\documentclass[11pt]{article}
\usepackage[a4paper,left=18mm,top=20mm,right=18mm,bottom=25mm]{geometry}

\usepackage{authblk}
\usepackage[numbers, sort&compress]{natbib}
\usepackage[margin=30pt,font=small,labelfont=bf]{caption}
\usepackage{amsmath,amsopn,amsthm,amssymb}

\usepackage{graphicx}
\usepackage{mathptmx} 
\usepackage[mathscr]{euscript}
\usepackage{natbib}
\usepackage{mathtools}
\usepackage{enumitem}
\usepackage{xspace}
\usepackage{scalerel}
\usepackage{bold-extra}

\usepackage[
 	colorlinks=true,
	linkcolor=blue
]{hyperref}

\newcommand{\ANC}{\ensuremath{\mathtt{ANC}}}
\newcommand{\LCA}{\ensuremath{\operatorname{LCA}}}
\newcommand{\lca}{\ensuremath{\operatorname{lca}}}

\newcommand{\lop}{\ensuremath{\operatorname{lop}}}
\newcommand{\lxt}{\ensuremath{\operatorname{lxt}}}
\renewcommand{\sf}{\ensuremath{\operatorname{sf}}}
\newcommand{\rev}{\ensuremath{\operatorname{rev}}}

\newcommand{\hasse}{\mathscr{H}}
\newcommand{\Hasse}{\mathscr{H}}
\newcommand{\DD}[2]{\ensuremath{\mathtt{D}_{#2}(#1)}} 
\newcommand{\desccl}[1]{\ensuremath{\mathfrak{D}_{#1}}}

\newcommand{\predecessorscl}[1]{\ensuremath{\mathfrak{A}_{#1}}}
\newcommand{\allpath}[3]{\ensuremath{\mathtt{B}_{#3}(#1, #2)}}
\newcommand{\allpathcl}[1]{\ensuremath{\mathfrak{I}_{#1}}}

\DeclareMathOperator{\parent}{parent}
\DeclareMathOperator{\indeg}{indeg}
\DeclareMathOperator{\outdeg}{outdeg}
\DeclareMathOperator{\cl}{cl}

\newcommand{\One}{\mathscr{I}}

\newcommand{\Ktt}{\ensuremath{K_{2,2}}}
\newcommand{\holju}{\textsc{holju}\xspace}

\DeclareMathOperator{\CC}{\mathtt{C}}
\DeclareMathOperator{\BB}{\mathtt{B}}

\newtheorem{theorem}{Theorem}[section]
\newtheorem{lemma}[theorem]{Lemma}
\newtheorem{corollary}[theorem]{Corollary}
\newtheorem{proposition}[theorem]{Proposition}

\newtheorem{definition}[theorem]{Definition}
\newtheorem{observation}[theorem]{Observation}
\newtheorem{remark}{Remark}

\newcommand{\REV}[1]{}

\providecommand{\keywords}[1]{\textbf{\textit{Keywords: }} #1}

\title{Global Least Common Ancestor (LCA) Networks\footnote{In memory of Andreas Dress, whose pioneering work and dedication to discrete mathematics and its interdisciplinary applications - among many other fields - have inspired generations of researchers. His contributions will undoubtedly continue to influence and guide our work for many years to come.}}

\author[1]{Anna Lindeberg}
\author[3,4]{Bruno J.\ Schmidt}
\author[2]{Manoj Changat}
\author[2]{Ameera Vaheeda Shanavas}
\author[3,4,5,6,7,8]{Peter F.\ Stadler}
\author[1,$\ddagger$]{Marc Hellmuth}

\affil[1]{Department of Mathematics, Faculty of Science,
  Stockholm University, SE-10691 Stockholm, Sweden} 

\affil[2]{Department of Futures Studies, University of Kerala, Karyavattom,
  695 581 Thiruvananthapuram, India}

\affil[3]{Max Planck Institute for Mathematics in the Sciences,
      Inselstra{\ss}e 22, D-04103 Leipzig, Germany}

\affil[4]{Bioinformatics Group, Department of Computer Science \&
    Interdisciplinary Center for Bioinformatics, Leipzig
      University, H{\"a}rtelstra{\ss}e 16–18,
      D-04107 Leipzig, Germany}

\affil[5]{Zuse School for Embedded and Composite Artificial
  Intelligence (SECAI) \& Center for Scalable Data Analytics and
  Artificial Intelligence (ScaDS.AI)
  Dresden-Leipzig \&
  German Centre for Integrative Biodiversity Research (iDiv)
  Halle-Jena-Leipzig \& Centre for Biotechnology and Biomedicine \&
  LIFE - Leipzig Research Centre for Civilization Diseases,
  Leipzig University, Leipzig, Germany}

\affil[6]{Department of Theoretical Chemistry, University
    of Vienna, W{\"a}hringerstra{\ss}e 17,
    A-1090 Wien, Austria}

\affil[7]{Facultad de Ciencias, Universidad National de
    Colombia, Bogot{\'a}, Colombia}

\affil[8]{Santa Fe Institute, 1399 Hyde Park Rd., Santa Fe, NM 87501, USA}

\affil[$\ddagger$]{Corresponding author: \textnormal{\texttt{marc.hellmuth@math.su.se}}	}

\date{\ }

\begin{document}
\sloppy

\maketitle

\abstract{ Directed acyclic graphs (DAGs) are fundamental structures used
  across many scientific fields. A key concept in DAGs is the least common
  ancestor (LCA), which plays a crucial role in understanding hierarchical
  relationships. Surprisingly little attention has been given to DAGs that
  admit a unique LCA for every subset of their vertices.  Here, we
  characterize such global lca-DAGs and provide multiple structural and
  combinatorial characterizations. We show that global lca-DAGs have a
  close connection to join semi-lattices and establish a connection to
  forbidden topological minors. In addition, we introduce a constructive
  approach to generating global lca-DAGs and demonstrate that they can be
  recognized in polynomial time. We investigate their relationship to
  clustering systems and other set systems derived from the underlying
  DAGs.
}

\smallskip
\noindent


\section{Introduction}

Directed acyclic graphs (DAGs) are widely used models in many fields of
application, ranging from algorithms and data structures
\cite{aho2006compilers,muchnick1997compiler,cormen2022introduction},
machine learning and artificial intelligence \cite{pearl2009causality,
  koller2009graphical}, to biology and bioinformatics
\cite{davidson2006regulatory, Huson:11,Nakhleh:05}, operations research and
decision making \cite{pinedo2016scheduling}, social and economic sciences
\cite{FH:14, ahuja1993networkflows}, and transportation and logistics
\cite{ahuja1993networkflows}. Each DAG $G = (V,E)$ naturally gives rise to
a partially ordered set $(V,\preceq)$, defined by $u \preceq v$ for
vertices $u,v\in V$ whenever there is a directed path from $v$ to $u$ in
$G$. In particular, each vertex $v \in V$ is naturally associated with the
set $\DD{v}{}\coloneqq\{u\in V(G)\mid u\preceq v\}$ of its descendant
vertices, i.e., those reachable from $v$. The set of leaf-descendants of
$v$ is then defined by $\CC(v)=\DD{v}{}\cap L(G)$ where $L(G)$ is the set
of all $\preceq$-minimal vertices{, i.e., the leaves} in $G$.  A least
common ancestor (LCA) of a subset $A \subseteq V(G)$ in a DAG $G$ is a
$\preceq$-minimal vertex $v$ in $G$ such that $A \subseteq \DD{v}{}$.

The existence of a \emph{unique} LCA for subsets $A \subseteq V$ plays a
central role in numerous applications. In evolutionary networks $G$, for
example, a unique LCA of $A$ represents a well-defined ancestral state from
which the taxa in $A$ have diverged, and is therefore crucial for the
interpretation of evolutionary relationships
\cite{Hellmuth:22q,HL:24,NW:05}. {Further applications arise in the
  identification of the most direct common cause in Bayesian networks
  \cite{pearl2009causality}, in object inheritance within programming
  languages \cite{AHBLN:89}, and in algorithms on lattices
  \cite{Jourdan1994,AHBLN:89}; see \cite{thesisJN:09} for an overview.}

{Much of the existing structural theory is concerned with
  \emph{leaf-based} LCA questions, which naturally arise in
    particular in phylogenetic applications \cite{Huson:11}. In this
    context, rooted trees and networks encode 
  hypothetical evolutionary histories, while empirical data are
  almost exclusively confined to extant species,
  i.e., the leaves. It is therefore natural
  that many structural results are formulated in terms of subsets
  $A\subseteq L(G)$ and the associated leaf-descendant clusters, LCAs, or
  related common-ancestor concepts
  \cite{Nakhleh:05,Huson:11,Hellmuth:22q,HL:24,HL:24a,FISCHER2010331,vanIersel2014}. For
  instance, \cite{S+24} studies conditions ensuring unique LCAs for leaf
  subsets of cardinality at most $k$ and relates them to $k$-weak
  hierarchies, and in \cite{HL:24a} classes of DAGs have been characterized
   in which every non-empty subset of leaves  admit a unique LCA. 
  }

{In contrast, the present work is concerned with the genuinely
  \emph{global} question: which DAGs admit a unique LCA for \emph{every}
  non-empty subset of vertices, not only for subsets of leaves?  This
  distinction is essential: leaf-based conditions do not in general
  determine the corresponding vertex-based LCA structure. Related
  computational problems have been studied for arbitrary vertex subsets, in
  particular for pairs of vertices: numerous results address the task of
  determining an element of $\LCA(A)$ when $A\subseteq V(G)$ consists of
  $|A|=2$ vertices
  \cite{BENDER200575,Grandoni:21,Mathialagan:22}. Generalizations such as
  deciding whether two vertices admit a unique LCA and computing the
  complete set of LCAs for a given pair are discussed in
  \cite{EMAN:07,Kowaluk:07,BEGKN:07}, and extensions to subsets of size
  $|A|\geq 3$ have also been considered
  \cite{thesisJN:09,YUSTER2008145,KLN:09}. Our main goal here is
  complementary to these algorithmic questions: we characterize those DAGs
  for which the LCA is uniquely defined for all non-empty subsets of
  $V(G)$.}

{This contribution is organized as follows. After introducing notation
  and preliminary facts, we derive in Section~\ref{sect:semilattice} the
  order-theoretic characterization via join-semilattices and prove that the
  global lca-property is equivalent to its pairwise
  version. Section~\ref{sec:Ktt} provides a characterization in terms of
  topological minors. In Section~\ref{sec:holju}, we give a constructive
  characterization showing that every global lca-DAG can be built
  recursively from a single-vertex graph. Section~\ref{sec:lxtG}
  reformulates the global lca-property in terms of the leaf-extended DAG
  $\lxt(G)$. In Section~\ref{sec:DAB-relations}, we characterize global
  lca-DAGs through closedness of descendant, ancestor, and intermediary set
  systems.  Section~\ref{sec:canonHasse} shows how a DAG can be
  reconstructed, up to shortcuts, from a canonical Hasse-diagram
  representation of an associated set system. We conclude in
  Section~\ref{sec:outlook} with a brief summary and open questions.}

\section{Notation and preliminary results}
\label{sec:basics}

\paragraph{Sets}
All sets considered here are assumed to be finite. The set $2^X$ denotes
the powerset of a set $X$.  A \emph{set system $\mathfrak{C}$ (on $X$)} is
a subset $\mathfrak{C} \subseteq 2^X$. A set system $\mathfrak{C}$ on $X$
is \emph{grounded} if $\{x\}\in \mathfrak{C} $ for all $x\in X$ and
$\emptyset\notin \mathfrak{C}$. A set system $\mathfrak{C}$ on $X$ is
\emph{rooted} if $X\in \mathfrak{C}$. A rooted and grounded set system is a
\emph{clustering system}.

Most of the set system considered here (e.g.\ those derived from DAGs)
will not contain the empty set.  Therefore, closedness is
naturally defined in terms of \emph{non-empty} intersections:
\begin{definition}
  \label{def:closed}
  A set system $\mathfrak{C}$ is closed if and only if $A,B\in
  \mathfrak{C}$ implies $A\cap B\in \mathfrak{C}$ or $A\cap B=\emptyset$.
\end{definition}
In the theory of convex structures \cite{vandeVel:93}, unlike here, set
systems are assumed to contain the empty set. This leads to a few subtle,
but rather trivial, differences in the definition of closure and closedness
that concern the role of the empty set. For instance, including the empty
set in $\mathfrak{C}$ relieves the need to treat $A\cap B=\emptyset$ as
a separate case in Definition~\ref{def:closed}.  We defer this point to at the
end of Section~\ref{sec:lxtG}, see also \cite{Changat:19a}.

A \emph{poset} $(\mathfrak{S},\sqsubseteq)$ is a set $\mathfrak{S}$  equipped with a 
partial order $ \sqsubseteq$. Two posets $(\mathfrak{S}_1,\sqsubseteq_1)$
and $(\mathfrak{S}_2,\sqsubseteq_2)$ are \emph{isomorphic} if there is a 
bijective map $\varphi \colon \mathfrak{S}_1 \to \mathfrak{S}_2$
that satisfies, for every $s,t\in \mathfrak{S}_1$, $s \sqsubseteq_1 t$
if and only if $\varphi(s) \sqsubseteq_2 \varphi(t)$.
A \emph{join-semilattice} is a poset
$(\mathfrak{S},\sqsubseteq)$ such that for all $s_1,s_2 \in \mathfrak{S}$
there is a least upper bound $s\coloneqq \sup_{\mathfrak{S}}(s_1,s_2)\in
\mathfrak{S}$, i.e., $s_1,s_2\sqsubseteq s$ and $s_1,s_2\sqsubseteq s'$
with $s'\in \mathfrak{S}$ implies $s \sqsubseteq s'$ \cite{Gratzer:07}.

\paragraph{DAGs}

A \emph{directed graph} $G=(V,E)$ is an ordered pair consisting of a non-empty set $V(G)\coloneqq V$ of
\emph{vertices} and a set $E(G)\coloneqq E \subseteq\left(V\times V\right)\setminus\{(v,v) \mid v\in
V\}$ of \emph{edges}.  
For directed graphs $G=(V_G,E_G)$ and $H=(V_H, E_H)$, an
\emph{isomorphism between $G$ and $H$} is a bijective map $\varphi\colon V_G\to V_H$ such that 
$(u,v)\in E_G$ if and only if $(\varphi(u),\varphi(v))\in E_H$. 
If such a map exists, then $G$ and $H$ are \emph{isomorphic}, in symbols $G\simeq H$.

A \emph{subgraph} of $G=(V,E)$ is a directed graph $H=(V',E')$ such that $V'\subseteq V$ and
$E'\subseteq E$. We occasionally denote this by $H\subseteq G$. If $E'$ comprises all edges $(u,v)$
in $E$ with $u,v\in V'$, then $H$ is the \emph{subgraph induced by $V'$}. In particular, $G-v$
denotes the subgraph of $G$ induced by $V\setminus\{v\}$, provided that $v\in V$ and $|V|>1$. A
\emph{directed $v_1v_n$-path} $P = (V,E)$ has vertex set $V = \{v_1,v_2,\ldots,v_n\}$ and the edges
in $E$ are precisely of the form $(v_i,v_{i+1})$, $i=1,2,\ldots, n-1$. If there is an edge
$(v_n,v_1)$ and a directed $v_1v_n$-path in $G$, then $G$ contains a \emph{cycle}. A directed graph
without cycles is a \emph{directed acyclic graph (DAG)}.

For a directed graph $G=(V,E)$, we define $\indeg_G(v)\coloneqq \left|\left\{u\in V \mid (u,v)\in
E\right\}\right|$ and $\outdeg_G(v)\coloneqq \left|\left\{u \in V\mid (v,u)\in E\right\}\right|$
for each $v\in V$ as the \emph{in-degree} respectively \emph{out-degree} of $v$ in $G$.

{Each DAG $G$ can be associated with}
a partial order $\preceq_G$ on the vertex set $V(G)$, defined by
$v\preceq_G w$ if and only if there is a directed $wv$-path. In
this case, we say that $w$ is an \emph{ancestor} of $v$ and $v$ is a \emph{descendant} of $w$. If $v\preceq_G w$
and $v\neq w$, we write $v\prec_G w$. If $(u,v) \in E(G)$, then $u$ is a \emph{parent} of $v$ and 
$v$ a \emph{child} of $u$. {We let $\parent_G(v)$ denote the set of all parents of $v$ in $G$.}
Two vertices $u,v\in V(G)$ are
\emph{$\preceq_{G}$-incomparable} if neither $u\preceq_G v$ nor $v\preceq_G u$ is true. We denote by
$L(G)\subseteq V(G)$ the $\preceq_G$-minimal vertices of $G$, and we call $x\in L(G)$ a \emph{leaf} of $G$.
If $L(G)=X$, then $G$ is a \emph{DAG on $X$}. 
A vertex $v\in V^0(G)\coloneqq V(G)\setminus L(G)$ is called an \emph{inner vertex}.
A vertex $v\in V(G)$ of a DAG $G$ that is $\preceq_G$-maximal is called a \emph{root}, and the set of roots of $G$ is denoted by $R(G)$. 
Note that $L(G)\neq \emptyset$ and $R(G)\neq \emptyset$ for all DAGs $G$ \cite{HL:24}.
A \emph{(rooted) network} $N$ is a DAG for which $|R(N)|=1$, i.e., $N$
has a unique root $\rho\in V(N)$. 

An edge $e=(u,w)$ in a DAG $G$ is a \emph{shortcut} if there is a vertex $v\in V(G)$ such that $w \prec_G v \prec_G u$, 
i.e., there is a directed $uw$-path that does not contain the edge $e$. 
A DAG without shortcuts is \emph{\textbf{s}hortcut-\textbf{f}ree ($\sf$)}. We denote with $\sf(G)$ the DAG obtained from $G$ by removing all
shortcuts.

For a DAG $G$ we define $\rev(G) \coloneqq (V(G), \{(v,u)\mid (u,v)\in E(G)\})$ as the \emph{reverse} of $G$.
Note that a directed graph $G$ is a DAG if and only if it admits a topological order $\lessdot_G$ defined on $V(G)$, i.e., 
the vertices of $G$ can be totally ordered such that $(x,y)\in E(G)$ implies that $x\lessdot_G y$ \cite{Kahn:62}.
It is easy to see $\lessdot_{\rev(G)}$ defined by  $x \lessdot_{\rev(G)} y$  if and only $y \lessdot_G x$
for all $x,y\in V(G)$ yields a  topological order on $\rev(G)$. Thus, we obtain

\begin{observation}\label{obs:rev-rev-G-is-G}
  Let $G$ be a DAG. Then $\rev(G)$ is a DAG and $\rev(\rev(G)) = G$.
\end{observation}

\paragraph{Clustering and Descendant Systems}
For every $v\in V(G)$ in a DAG $G$, the set of its descendant leaves
\begin{equation*}
  \CC_G(v)\coloneqq\{ x\in L(G)\mid x \preceq_G v\}
\end{equation*}
is a \emph{cluster} of $G$.  We write
$\mathfrak{C}_G\coloneqq\{\CC_G(v)\mid v\in V(G)\}$ for the set of all
clusters in $G$.  {Note that $\mathfrak{C}_G$ is grounded for all DAGs
  $G$. However, if a DAG $G$ has several roots, $\mathfrak{C}_G$ may not be
  rooted and, therefore, not a clustering system. For later reference, we
  provide the following results which shows, among others, that the set of
  clusters $\mathfrak{C}_G$ is invariant under shortcut-removal from $G$.}
\begin{lemma}[{\cite[L.~2.5]{HL:24}}]\label{lem:shortcutfree}
  Let $G = (V,E)$ be a DAG. Then, for all $u,v\in V$, it holds that
  $u\prec_G v$ if and only if $u\prec_{\sf(G)} v$ and, for all $v\in V$, it
  holds that $\CC_G(v)= \CC_{\sf(G)}(v)$.  In particular, $\sf(G)$ is a
  shortcut-free DAG.
\end{lemma}

\begin{lemma}[{\cite[L.~1]{S+24}}]\label{lem:prec-subset}
  For all DAGs $G$ and all $u,v\in V(G)$ it holds that $u\preceq_G v$
  implies $\CC_G(u)\subseteq\CC_G(v)$.
\end{lemma}

\begin{figure}
  \centering
  \includegraphics[width=0.9\textwidth]{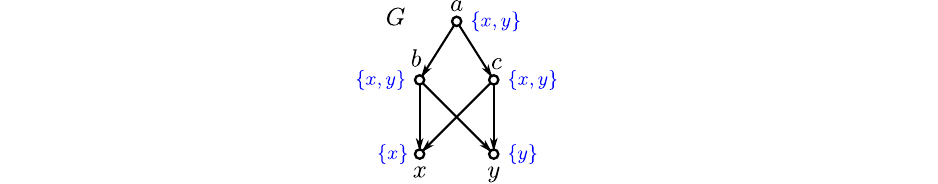}
  \caption{The clustering system $\mathfrak{C}_G =
    \{\{x\},\{y\},\{x,y\}\}$, indicated in blue, of $G$ is closed but $G$
    does not have the lca-property as $|\LCA_G(\{x,y\})|>1$.  Note that $G$
    contains a $\Ktt$-minor without ${\mathsf{X}}$- or
    ${\mathsf{X'}}$-subdivision, see Section~\ref{sec:Ktt} for more
    details. Although $\mathfrak{C}_G$ is closed, the descendant system
    $\mathfrak{D}_G$ of $G$ is not closed because $\DD{b}{G} \cap \DD{c}{G}
    = \{b,x,y\} \cap \{c,x,y\} = \{x,y\}$ but $\{x,y\} \notin
    \mathfrak{D}_G$.}
  \label{fig:exmpl-lca-closed}
\end{figure}

There are DAGs $G$ containing $\preceq_G$-incomparable vertices $u$ and $v$
for which $\CC_G(u)\subseteq\CC_G(v)$ is satisfied. For example, the
vertices $b$ and $c$ of the DAG $G$ in Figure~\ref{fig:exmpl-lca-closed}
are $\preceq_G$-incomparable but satisfy $\CC_G(b)=\CC_G(c)$.  The
following property ensures that there is a directed path between two
vertices $v$ and $u$ precisely if their clusters are comparable with
respect to inclusion.

\begin{definition}[{\cite{Hellmuth:22q}}]
  A DAG $G$ has the path-cluster-comparability (PCC) property if it
  satisfies, for all $u, v \in V(G)$: $u$ and $v$ are
  $\preceq_G$-comparable if and only if $\CC_G(u) \subseteq \CC_G(v)$ or
  $\CC_G(v) \subseteq \CC_G(u)$.
\end{definition}

For a DAG $G$ and a vertex $v\in V(G)$, let 
\begin{equation*}
	\DD{v}{G} \coloneqq \{ u\in V(G) \mid u \preceq_G v \}
\end{equation*}	
 be the set of all \emph{descendants} of $v$ in
$G$. Moreover, {the \emph{descendant system}} $\mathfrak{D}_G \coloneqq \{\DD{v}{G} \mid v\in V(G)\}$
{consists} of all descendant sets in $G$.  In general $\mathfrak{D}_G$
is not grounded as, for all non-leaf vertices $v\in V(G)\setminus L(G)$, we
have $\{v\}\notin \mathfrak{D}_G$ since $\{v\} \neq \DD{u}{G}$ for any
$u\in V(G)$.  Note also that for every vertex $v\in V(G)$ it holds that
$\CC_G(v)\subseteq\DD{v}{G}$, with equality if and only if $v$ is a
leaf. We note in passing that a construction resembling $\mathfrak{D}_G$
appears in \cite[Sec. 6]{Halmos:1974}, albeit in a quite different context.
Moreover, the descendant set $\mathfrak{D}_G$ of $G$ coincides with the
collection of all so-called principal lower sets of the poset
$(V(G),\preceq_G)$, see e.g. \cite{Gierz:80}.  The latter has, to the best
of our knowledge, so far not been studied in the context of DAGs and least
common ancestors. To keep this contribution self-consistent, we include
proofs of some of the properties of $\mathfrak{D}_G$ that also appear as
Exercise 1.16(ii) of \cite{Gierz:80}.

\begin{lemma}\label{lem:desc-subset-iff-comparable}
  Let $G$ be a DAG and $u,v\in V(G)$. Then, $\DD{u}{G} \neq \DD{v}{G}$ if
  and only if $u\neq v$.  Moreover, $u\preceq_G v$ if and only if
  $\DD{u}{G}\subseteq \DD{v}{G}$. In particular, $\DD{u}{G}\subsetneq
  \DD{v}{G}$ implies that $v\notin \DD{u}{G}$.  Consequently, the poset
  $(\mathfrak{D}_G,\subseteq)$ is isomorphic to the poset
  $(V(G),\preceq_G)$. 
\end{lemma}
\begin{proof}
  Let $G$ be a DAG and $u,v\in V(G)$. Clearly, $\DD{u}{G} \neq \DD{v}{G}$
  implies $u\neq v$. Assume that $u\neq v$. One of these cases must occur:
  (i) $u\prec_G v$, (ii) $v\prec_G u$ or (iii) $u$ and $v$ are
  $\preceq_G$-incomparable. Since $G$ is DAG, in Case (i) we have $v\not
  \preceq_G u$ and in Case (ii) we have $u\not \preceq_G v$. Thus, in Case
  (i), we have $v\in \DD{v}{G}$ but $v\notin \DD{u}{G}$ and, in Case (ii),
  we have $u\in \DD{u}{G}$ but $u\notin \DD{v}{G}$. In Case (iii) we have
  $u\in \DD{u}{G}$ but $u\notin \DD{v}{G}$. Hence, in all cases $\DD{v}{G}
  \neq \DD{u}{G}$.
	
  Assume now that $u\preceq_G v$. In this case, all $w\preceq_G u$
  trivially satisfy $w\preceq_G v$ and $\DD{u}{G}\subseteq \DD{v}{G}$
  follows.  Suppose now that $\DD{u}{G}\subseteq \DD{v}{G}$.  Since, by
  definition, $u\in \DD{u}{G}$ it follows that $u\in \DD{v}{G}$.
  Consequently, $u\preceq_G v$ must hold.  Finally, assume that
  $\DD{u}{G}\subsetneq \DD{v}{G}$. Since $u\in \DD{u}{G}$, we have $u\in
  \DD{v}{G}$. Hence, $u\prec_G v$. This together with the fact that $G$ is
  a DAG implies that $v\not \prec_G u$ and thus, $v\notin \DD{u}{G}$.

  Finally, the previous arguments also show that the map $\varphi:
  V(G)\to\mathfrak{D}_G$ defined by $v\mapsto\DD{v}{G}$ is a bijection that
  satisfy $u\preceq_G v$ if and only if $\DD{u}{G}\subseteq\DD{v}{G}$. In
  other words, $\varphi$ is an order-preserving bijection and thus an
  isomorphism of the posets $(V(G),\preceq_G)$ and
  $(\mathfrak{D}_G,\subseteq)$.
\end{proof}

There is a close connection between the {descendant system}  $\mathfrak{D}_G$ {and the set of clusters} $\mathfrak{C}_G$. Note that
$\mathfrak{D}_G$ is a set system on $V(G)$, while $\mathfrak{C}_G$ is a set
system on $L(G)\subseteq V(G)$. In particular, $x\in L(G)$ if and only if
$\{x\} = \DD{x}{G} \in \mathfrak{D}_G$. Moreover, one easily verifies that
$\CC_G(v) = \DD{v}{G} \cap L(G)$ for every $v\in V(G)$ and, therefore, that
$C \in \mathfrak{C}_G$ if and only if $C= D\cap L(G)$ for some $D \in
\mathfrak{D}_G$.  Therefore, we obtain
\begin{observation}\label{obs:recoverC}
  For every DAG $G$ it holds that $\mathfrak{C}_G=\{D\cap
  L(G)\,:\,D\in\mathfrak{D}_G\} $.
\end{observation}
By Observation~\ref{obs:recoverC}, it is straightforward of obtain the
{set of clusters} $\mathfrak{C}_G$ from the descendant system
$\desccl{G}$. As we shall see later in Section~\ref{sec:DAB-relations}, it
is also possible to construct $\desccl{G}$ from the {set of clusters} of
so-called leaf-extended versions $\lxt(G)$ of a DAG $G$, although the
procedure is less simple (cf.\ Lemma~\ref{lem:D-and-C_N*}).

{The descendant system of a DAG $G$ contains} much more information
than the corresponding {set of clusters} {of $G$}. Still,
{descendant systems} neglect some of the structure of $G$. In
particular, we have $\mathfrak{D}_G = \mathfrak{D}_{G'}$ whenever $G'$ is
obtained from $G$ by adding shortcuts
(cf.\ Lemma~\ref{lem:shortcutfree}). Nevertheless, every DAG $G$ is
uniquely determined by $\mathfrak{D}_G$ up to its shortcuts, as shown next.
\begin{lemma}\label{lem:sf(G)=hasseDG}
  For a DAG $G$ we have $V(G) = \cup_{D\in \mathfrak{D}_G} D$.  Let
  $E'\subseteq E(G)$ denote the set of all edges in $G$ that are not
  shortcuts.  Then, for any two vertices $u,v\in V(G)$, it holds that
  $(u,v)\in E'$ if and only if $v\in \DD{u}{G}$, and for all $w\in
  \DD{u}{G}\setminus\{u,v\}$, we have $v\notin \DD{w}{G}$.
\end{lemma}
\begin{proof}
  Let $G$ be a DAG.  Since $v\in \DD{v}{G}$ and $\DD{v}{G} \in
  \mathfrak{D}_G$ for all $v\in V(G)$ it follows that $V(G) = \cup_{v\in
    V(G)} \DD{v}{G} = \cup_{D\in \mathfrak{D}_G} D $.  Let $E'\subseteq
  E(G)$ denote the set of all edges in $G$ that are not shortcuts. Let
  $(u,v)\in E'$. Clearly, $v\in \DD{u}{G}$ holds. Let $w\in \DD{u}{G}$ with
  $w\notin\{u,v\}$.  Since $(u,v)$ is not a shortcut in $G$, it cannot hold
  that $v\prec_G w\prec_G u$.  Hence, $v\notin \DD{w}{G}$.  Suppose now
  that $u,v\in V(G)$ are distinct vertices for which it holds that $v\in
  \DD{u}{G}$ and for all $w\in \DD{u}{G}\setminus\{u,v\}$ we have $v\notin
  \DD{w}{G}$. Hence, $v\prec_G u$ and no descendant $w\notin\{u,v\}$ of $u$
  is an ancestor of $v$ in $G$. Hence, there is no vertex $w\in V(G)$ that
  satisfies $v\prec_G w\prec_G u$. This and $v\prec_G u$ now implies that
  $(u,v)\in E(G)$ and, in particular, that $(u,v)$ is not a shortcut.
\end{proof}
    
\paragraph{Least Common Ancestors (LCAs)}

For a given DAG $G$ and a non-empty subset $A\subseteq V(G)$, a vertex $v\in
V(G)$ is a \emph{common ancestor of $A$} if $v$ is an ancestor of every
vertex in $A$.  We denote with $\ANC_G(A)$ the set of all common ancestors
of $A$ in $G$.  Moreover, $v$ is a \emph{least common ancestor} (LCA) of
$A$ if $v$ is a $\preceq_G$-minimal common ancestor of $A$.  The set
$\LCA_G(A)$ comprises all LCAs of $A$ in $G$. In a network $N$, the unique
root is a common ancestor for all non-empty $A\subseteq V(N)$ and,
therefore, $\LCA_N(A)\neq\emptyset$ for all non-empty $A\subseteq V(N)$.  
In this contribution, we 
{will encounter several situations in which}
there is a unique LCA of certain subsets of the vertices, i.e.\ when $|\LCA_G(A)|=1$
for {some non-empty subset} $A\subseteq V(G)$.  For simplicity, we will write $\lca_G(A)=v$ in case
that $\LCA_G(A)=\{v\}$ and say that \emph{$\lca_G(A)$ is well-defined};
otherwise, we leave $\lca_G(A)$ \emph{undefined}. {Moreover, we will write $\lca_G(u,v)$ instead of $\lca_G(\{u,v\})$ for $u,v\in V(G)$.}

{The following definition specifies DAGs where the LCA is unique for all non-empty (size-2) subsets of the \emph{leaf set}. They have been explicitly studied in \cite{Hellmuth:22q,HL:24a,S+24}. }
\begin{definition}
  Let $G = (V,E)$ be a DAG on $X$. Then, 
  \begin{itemize}[noitemsep]
  \item $G$ has the \emph{lca-property} if $|\LCA_G(A)|=1$ for all
    non-empty subsets $A\subseteq X$.
  \item $G$ has \emph{pairwise-lca-property} if $|\LCA_G(A)|=1$ for all
    subsets $A\subseteq X$ of size $|A|=2$.
  \end{itemize}
\end{definition}
\noindent Networks that have the lca-property are also called
\emph{lca-networks}.  Moreover, $G$ is \emph{lca-relevant} if, for all
$v\in V(G)$, there is some subset $A\subseteq L(G)$ with $v=\lca_G(A)$. %
We emphasize that a DAG can be lca-relevant without having the
lca-property, and \textit{vice versa}. For examples and details on the
relationship between these two properties we refer to \cite{HL:24,HL:24a}.
The following result generalizes Lemma~39 in \cite{Hellmuth:22q}, extending
the result from networks to general DAGs.

\begin{lemma}\label{lem:lca->closed} 
  If a DAG $G$ has the lca-property, then $\mathfrak{C}_G$ is a closed
  clustering system.
\end{lemma}
\begin{proof}
  Let $G$ be a DAG on $X$ with the lca-property. Clearly,
  $\emptyset\notin\mathfrak{C}_G$ and $\{x\}\in\mathfrak{C}_G$ for each
  $x\in X$ and, thus, $\mathfrak{C}_G$ is grounded. To show that
  $\mathfrak{C}_G$ is a clustering system it thus suffices to show that
  $X\in\mathfrak{C}_G$, i.e., that $\mathfrak{C}_G$ is rooted. Since $G$
  has the lca-property, $v=\lca_G(X)$ must be well-defined. By definition
  of least common ancestors, $X\subseteq \CC_G(v)$ which, together with
  $\CC_G(v)\subseteq X$, implies that $\CC_G(v) = X$. Hence, $X\in
  \mathfrak{C}_G$. In summary, $\mathfrak{C}_G$ is a clustering system.

  It remains to show that $\mathfrak{C}_G$ is closed. Let $\CC_G(u),
  \CC_G(v) \in \mathfrak{C}_G$ for some $u,v\in V(G)$. If $\CC_G(u)
  \subseteq \CC_G(v)$, $\CC_G(v) \subseteq \CC_G(u)$ or $\CC_G(u) \cap
  \CC_G(v)=\emptyset$, there is nothing to show. Hence, assume that
  $\CC_G(u)\cap\CC_G(v)\notin\{\CC_G(u),\CC_G(v),\emptyset\}$ and put
  $A\coloneqq\CC_G(u)\cap \CC_G(v)\ne\emptyset$. Since $G$ has the
  lca-property, there is $w\in V(G)$ such that $w=\lca_G(A)$, and thus
  $A\subseteq \CC_G(w)$. The contraposition of Lemma~\ref{lem:prec-subset}
  shows that $u$ and $v$ are two $\preceq_G$-incomparable common ancestors
  of $A$. Since $w$ is the unique $\preceq_G$-minimal common ancestor of
  $A$, we have $w\preceq_G u$ and $w\preceq_G v$, which together with
  Lemma~\ref{lem:prec-subset} implies $\CC_G(w)\subseteq \CC_G(u)$ and
  $\CC_G(w)\subseteq \CC_G(v)$. Therefore $\CC_G(w)\subseteq A$.  Hence
  $A=\CC_G(w)\in\mathfrak{C}_G$ and thus, $\mathfrak{C}_G$ is closed.
\end{proof}

The converse of Lemma~\ref{lem:lca->closed} is not true. A counter-example
is shown in Figure~\ref{fig:exmpl-lca-closed}.

The next rather simple but technical result will be helpful to simplify
several arguments throughout this paper. 
\begin{lemma}\label{lem:common-anc-above-lca}
  Let $G$ be a DAG such that, for some $x,y\in V(G)$, $\lca_G(x,y)$ is
  well-defined.  If $u,v\in V(G)$ are $\preceq_G$-incomparable vertices of
  $G$ such that both $u$ and $v$ are common ancestors of $\{x,y\}$, then
  $\lca_G(x,y)\prec_G u$ and $\lca_G(x,y)\prec_G v$.
\end{lemma}
\begin{proof}
  By definition of least common ancestors, $\lca_G(x,y)\preceq_G v$ and
  $\lca_G(x,y)\preceq_G u$ holds. Therefore, if $\lca_G(x,y)\in\{u,v\}$,
  then $u$ and $v$ are $\preceq_G$-comparable. Since $u$ and $v$
  are $\preceq_G$-incomparable, it thus follows that $\lca_G(x,y)\prec_G v$
  and $\lca_G(x,y)\prec_G u$.
\end{proof}

For a poset $(\mathfrak{Q},\sqsubseteq)$, the \emph{Hasse diagram}
$\Hasse(\mathfrak{Q},\sqsubseteq)$ is the DAG with vertex set
$\mathfrak{Q}$ and directed edges from $A\in\mathfrak{Q}$ to
$B\in\mathfrak{Q}$ if (i) $B\sqsubseteq A$ and $A\neq B$ and (ii) there is
no {$C\in\mathfrak{Q}$} with $B\sqsubseteq C\sqsubseteq A$ and $C\neq A,B$.

A DAG $G=(V,E)$ is \emph{regular} if the map $\varphi\colon V\to
V(\Hasse(\mathfrak{C}_G,\subseteq))$ defined by $v\mapsto \CC_G(v)$ is an
isomorphism between $G$ and $\Hasse(\mathfrak{C}_G,\subseteq)$ where
$\subseteq$ denotes the usual subset-relation.  The next result summarizes
some of the results from \cite[Sec.\ 4]{HL:24}.
\begin{theorem}\label{thm:char-lcaRel}
  The following statements are equivalent for every DAG $G$. 
  \begin{enumerate}[noitemsep,label=(\arabic*)]
  \item $G$ is lca-relevant
  \item {$v=\lca_G(\CC_G(v))$ for  all $v\in V(G)$ (known as the \emph{strong-(CL) property})}
  \item $G$ satisfies (PCC) and $\CC_G(u)\neq \CC_G(v)$ for all distinct
    $u,v\in V(G)$.  \label{eq:PCC-neqCluster}
  \item $\CC_G(u)\subseteq \CC_G(v)$ if and only if $u\preceq_G v$ for all
    $u,v\in V(G)$. \label{eq:stronger-than-PCC}
  \item $G$ is isomorphic to the regular DAG
    $\Hasse(\mathfrak{C}_G,\subseteq)$ to which {$k\geq 0$} shortcuts have
    been added. \label{eq:sf-reg}
  \end{enumerate}
\end{theorem}

\section{Basic properties of global lca-networks and a first characterization}
\label{sect:semilattice}

The lca-property can be strengthened as follows. 

\begin{definition}\label{def:globallca}
  Let $G = (V,E)$ be a DAG on $X$. Then,
  \begin{itemize}[noitemsep]
  \item $G$ has the \emph{global lca-property} if $|\LCA_G(A)|=1$ for
    all non-empty subsets $A\subseteq V$. 
  \item $G$ has the \emph{global pairwise-lca-property} if
    $|\LCA_G(A)|=1$ for all subsets $A\subseteq V$ of size $|A|=2$.
  \end{itemize}
  Moreover, DAGs that have the global lca-property are also called
  \emph{global lca-DAGs}.
\end{definition}

The next result will be used throughout the paper without explicit mention
each time.
\begin{observation}
 Every global lca-network is an lca-network.
\end{observation}

Even more interesting is the following\
\begin{lemma}\label{lem:globallca=>network}
  If a DAG $G$ on $X$ has the global lca-property then $G$ is a network
  and $\mathfrak{C}_G$ is a closed clustering system on $X$.
\end{lemma}
\begin{proof}
  Let $G$ be a DAG on $X$ with global lca-property.  Suppose that $G$
  is not a network and thus, $|R(G)|>1$.  By definition,
  $\LCA_G(R(G))=\emptyset$. Hence, $G$ does not have the
  global lca-property; a contradiction.  Hence, $G$ must be a network.
  By \cite[L.~14]{Hellmuth:22q}, for any network $G$, the set system
  $\mathfrak{C}_G$ is a clustering system. Moreover, since $G$ has, in
  particular, the lca-property,  Lemma~\ref{lem:lca->closed} implies
  that $\mathfrak{C}_{G}$ is closed.
\end{proof}

We note that, similar to the ideas in \cite{HL:24}, we may define the
global LCA-property as a property of a DAG $G$ to ensure that for all
subsets (possibly of a given size) we have $\LCA(A)\neq \emptyset$ for all
$A\subseteq V(G)$, i.e., a not necessarily unique LCA for all subsets
$A\subseteq V(G)$ exists. However, by similar arguments that lead to
Lemma~\ref{lem:globallca=>network}, each DAG $G$ with the global 
LCA-property must be a network with a unique root $\rho_G$.  In this case,
however, $\CC_G(\rho_G) = L(G)$, which implies that $\rho_G$ serves as an
ancestor of all subsets $A$.  Consequently, $\LCA(A)\neq \emptyset$ for all
$A\subseteq V(G)$ is trivially satisfied if and only if $G$ is a network.
Thus, we focus here on the more interesting case of the global
lca-property as in Definition~\ref{def:globallca}.

\begin{remark}
  Since global lca-DAGs must be networks, we call such DAGs also
  \emph{global lca-networks}.
\end{remark}

The converse of Lemma~\ref{lem:globallca=>network} is not true.  In
particular, Figure~\ref{fig:exmpl-lca-closed} shows a DAG $G$, for which
$\mathfrak{C}_G$ is a closed clustering system but where $G$ does not have
the lca-property and thus, in particular, not the
global lca-property. As we shall see later in
Section~\ref{sec:DAB-relations}, closedness of $\mathfrak{D}_G$ instead
characterizes global lca-networks. However, in case that $G$
satisfies (PCC), closedness of its clustering system characterizes 
global lca-network as shown in the next result, which
generalizes Theorem 7 in \cite{Hellmuth:22q}. 
\begin{proposition}\label{prop:PPC+Closed=>globalLCA}
Let $G$ be a DAG satisfying (PCC). If $\mathfrak{C}_{G}$ is a clustering
system, then $G$ is a network. Moreover, $G$ is a global lca-network if and
only if $\mathfrak{C}_{G}$ is a closed clustering system.
\end{proposition}
\begin{proof}
  Let $G$ be a DAG satisfying (PCC).  Suppose first that
  $\mathfrak{C}_{G}$ is a clustering system.  Let $v\in V(G)$ be a vertex
  such that $\CC_{G}(v)=X$, which exists since $\mathfrak{C}_{G}$ is
  rooted.  Clearly, there is some root $r\in R(G)$ such that $v\preceq_G
  r$.  By Lemma~\ref{lem:prec-subset}, this vertex $r$ must satisfy $X
  =\CC_G(v)\subseteq \CC_G(r) \subseteq X$ and thus, $\CC_G(r) = X$.  Since
  $r$ is a root in $G$, there is no vertex $u$ satisfying $r\prec_G u$.
  This together with (PCC) and the fact that $\CC_{G}(u)\subseteq X$ for
  all vertices $u$ in $G$ implies that $u\preceq_{G} r$ for all vertices
  $u$ in $G$.  Hence $r$ is the unique root in $G$ and $G$ is a network.
  
  Suppose now that $\mathfrak{C}_{G}$ is a closed clustering system.  By
  the latter arguments, $G$ is a network and thus, $\LCA_G(A)\neq\emptyset$
  for all $A\subseteq V(G)$.  Assume, for contradiction, that there is a
  non-empty set $A\subseteq V(G)$ with $|\LCA_G(A)|>1$. Let
  $u,v\in\LCA_G(A)$ be two distinct vertices. Then $u,v$ must be
  $\preceq_G$-incomparable.  {Let $a\in A$. Since $u,v\in\LCA_G(A)$, it
    follows that $a\preceq_G u$ and $a\preceq_G v$. Lemma
    \ref{lem:prec-subset} implies that $\CC_G(a)\subseteq \CC_G(u)$ and
    $\CC_G(a)\subseteq \CC_G(v)$.  Hence, $C \coloneqq \CC_G(u)\cap
    \CC_G(v)\neq \emptyset$. Since $\mathfrak{C}_{G}$ is closed, $C\in
    \mathfrak{C}_{G}$. Thus, there is a vertex $w\in V(G)$ such that
    $\CC_G(w)=C$. Note that, since $G$ satisfies (PCC), all vertices
    $w',w''$ with $\CC_G(w')=\CC_G(w'')=C$ must be $\preceq_G$-comparable.
    Thus, we can assume without loss of generality that $w$ is the unique
    $\preceq_G$-maximal vertex satisfying $\CC_G(w)=C$.  Note that neither
    $\CC_G(u)\subseteq \CC_G(v)$ nor $\CC_G(v)\subseteq \CC_G(u)$ is
    possible because $G$ satisfies (PCC) and $u,v$ are
    $\preceq_G$-incomparable.  This and $\CC_G(u)\cap \CC_G(v) = \CC_G(w)$
    implies that neither $\CC_G(u)=\CC_G(w)$ nor $\CC_G(v)=\CC_G(w)$ can
    hold.  Hence, $\CC_G(w)\subsetneq \CC_G(u)$ and $\CC_G(w)\subsetneq
    \CC_G(v)$. In particular, $w\neq u$ and $w\neq v$.  The latter two
    arguments together with (PCC) and Lemma~\ref{lem:prec-subset} implies
    that $w\prec_G u$ and $w\prec_G v$.  Now let $b\in A$, where $b=a$ may
    be possible.  Since $u,v\in\LCA_G(A)$, it follows that $b\preceq_G u$
    and $b\preceq_G v$. Lemma \ref{lem:prec-subset} implies that
    $\CC_G(b)\subseteq \CC_G(u)$ and $\CC_G(b)\subseteq \CC_G(v)$.
    Therefore, $\CC_G(b)\subseteq \CC_G(u)\cap \CC_G(v) = \CC_G(w)$.  Since
    $G$ satisfies (PCC), $b$ and $w$ must be $\preceq_G$-comparable.  If $w
    \prec_G b$, then Lemma \ref{lem:prec-subset} implies that
    $\CC_G(w)\subseteq \CC_G(b)$ and thus, $\CC_G(w) = \CC_G(b)$; a
    contradiction to $w$ being the $\preceq_G$-maximal vertex satisfying
    $\CC_G(w)=C$.  Hence, $b \preceq_G w$ must hold. As $b\in A$ was chosen
    arbitrarily, it follows that for all $b\in A$ it holds that $b
    \preceq_G w$.  Hence, $w$ is a common ancestor of the elements in $A$.
    However, since $w\prec_G u$ and $w\prec_G v$, we have
    $u,v\notin\LCA_G(A)$; a contradiction.  In summary, $|\LCA_G(A)|=1$ for
    all non-empty subsets $A\subseteq V(G)$.  Therefore, $G$ is a global
    lca-network.}  Finally, if $G$ is a global lca-network, then
  Lemma~\ref{lem:globallca=>network} implies that $\mathfrak{C}_{G}$ is a
  closed clustering system.
\end{proof}
    
We now provide a first characterization of global lca-networks.  To this
end, we will adapt ideas established for semilattices
(cf.\ e.g.\ \cite[page 4]{Gratzer:07}).
    
\begin{theorem}   \label{thm:global-k2}
  The following statements are equivalent for every DAG $G$.
  \begin{enumerate}[label=(\arabic*)]
  \item $G$ has the global pairwise-lca-property
  \item  $G$ has the global lca-property.
  \end{enumerate}
\end{theorem}
\begin{proof}
  Clearly, Condition (2) implies (1).  We show now that Condition (1)
  implies (2). To this end, suppose that $G$ is a DAG that has the global
  pairwise-lca-property. Clearly, $\lca_{G}(\{x\})=x$ holds for all $x\in
  V(G)$.  Consider now $W=\{a_1,a_2,a_3\} \subseteq V(G)$. Hence,
  $\lca_{G}(a_1,a_2) = v$ and $\lca_{G}(a_3,v) = w$ exist. We claim that $w
  =\lca_{G}(\{a_1,a_2,a_3\})$. To see this, first observe that, by choice
  of $v$ and $w$, we have $a_1,a_2\preceq_{G} v$ respectively
  $v,a_3\preceq_{G} w$ which further implies that $a_1,a_2,a_3\preceq_{G}
  w$ and thus, $w$ is a common ancestor of $W$. Now, let $u$ be an
  arbitrary common ancestor of $a_1,a_2$ and $a_3$. Thus,
  $a_1,a_2,a_3\preceq_{G} u$ and therefore $u$ is, in particular, a common
  ancestor of $a_1$ and $a_2$. By definition of least common ancestors,
  $\lca_G(a_1,a_2)=v\preceq_{G} u$ holds.  {Since $v\preceq_{G} u$ and
    $a_3\preceq_{G} u$, the definition of least common ancestors
    furthermore implies that $w = \lca_{G}(a_3,v)\preceq_{G} u$.}  As these
  arguments hold for any common ancestor $u$ of $a_1,a_2,a_3$, we can
  conclude that $w = \lca_{G}(a_3,v) = \lca_{G}(a_1,a_2,a_3)$. Applying the
  latter arguments inductively, one arrives at $\lca_{G}(W) =
  \lca_{G}(\dots(\lca_{G}(a_1,a_2),a_3), \dots,a_\ell)$ for all non-empty
  sets $W = \{a_1,\dots, a_\ell\} \subseteq V(G)$. Hence, $G$ has the
  global lca-property.
\end{proof}

{As a consequence of Theorem~\ref{thm:global-k2}, it
  suffices to check the existence of a unique least common ancestor for all
  \emph{pairs} of vertices to conclude that a DAG has the global
  lca-property. Since Kowaluk and Lingas \cite[Thm.~3]{Kowaluk:07} have
  shown that verifying the existence of a unique LCA for all pairs of
  vertices of a given DAG can be achieved in polynomial time, we
  immediately obtain the following result.}

\begin{corollary}\label{cor:polytime}
  It can be verified in polynomial time if a DAG $G$
  has the global lca-property.
\end{corollary}

The remaining three results of this section formalize the, intuitively
clear, connections between join-semilattices and global lca-networks:
\begin{proposition}\label{prop:joinSL-hasse}
  For every poset $(\mathfrak{Q},\sqsubseteq)$ the following statements
  are equivalent.
  \begin{enumerate}[label=(\arabic*)]
  \item The Hasse diagram  $\Hasse(\mathfrak{Q},\sqsubseteq)$ has
    the global pairwise-lca-property.
  \item The Hasse diagram $\Hasse(\mathfrak{Q},\sqsubseteq)$ has the
    global lca-property.
  \item $(\mathfrak{Q},\sqsubseteq)$ is a join-semilattice.
  \end{enumerate}
\end{proposition}
\begin{proof}
  By Theorem~\ref{thm:global-k2}, Statement (1) and (2) are equivalent.
  Suppose that the Hasse diagram $\Hasse\coloneqq
  \Hasse(\mathfrak{Q},\sqsubseteq)$ has the global
  pairwise-lca-property.  Note that $X \preceq_{\Hasse} Y$ if and only
  if $X\sqsubseteq Y$ for all $X,Y \in \mathfrak{Q}$.  Moreover, since the
  Hasse diagram $\Hasse$ has the global pairwise-lca-property, for all
  $X,Y \in \mathfrak{Q}$, $Z_{XY}\coloneqq \lca_{\Hasse}(X,Y)$ is
  well-defined.  By definition, $Z_{XY} = \sup_{\mathfrak{Q}}(X,Y)$ and
  thus, $(\mathfrak{Q},\sqsubseteq)$ is a join-semilattice. Therefore, (1)
  implies (3).  Suppose that (3) holds. Thus, for $X_1,X_2 \in
  \mathfrak{Q}$ there is a unique least upper bound
  $\sup_{\mathfrak{Q}}(X_1,X_2)\in \mathfrak{Q}$. It now immediately
  follows that $\sup_{\mathfrak{Q}}(X_1,X_2)$ is a vertex in $\Hasse$.
  Moreover, by definition of the Hasse diagram, $\lca_{\Hasse}(X_1,X_2) =
  \sup_{\mathfrak{Q}}(X_1,X_2)$.  Hence, $\Hasse$ has the global
  pairwise-lca-property.
\end{proof}
 
\begin{corollary}\label{cor:joinSL-hasse}
  If a set system $\mathfrak{Q}$ is grounded and
  $(\mathfrak{Q},\subseteq)$ is a join-semilattice, then $\mathfrak{Q}$
  is a clustering system and there is a DAG $G$ with $\mathfrak{C}_G =
  \mathfrak{Q}$ and with global pairwise-lca-property.  In this case,
  $G$ is a network.
\end{corollary}
\begin{proof}
  Assume that $\mathfrak{Q}\subseteq 2^X$ is a grounded set system on $X$
  and $(\mathfrak{Q}, \subseteq)$ is a join-semilattice.  Since
  $\mathfrak{Q}$ is grounded, we have $\{x\}\in \mathfrak{Q}$ for all $x\in
  X$. Consider the graph $G \coloneqq \Hasse(\mathfrak{Q}, \subseteq)$,
  obtained by replacing each singleton vertex $\{x\}$ in
  $\Hasse(\mathfrak{Q},\subseteq)$ by the element $x$.  Since
  $L(\Hasse(\mathfrak{Q},\subseteq))=\{\{x\}\mid x\in X\}$ follows from
  $\mathfrak{Q}$ being grounded (in particular,
  $\emptyset\notin\mathfrak{Q}$) we have, together with the definition of
  Hasse diagrams, $L(G)=X$.  Moreover, the construction of $G$ ensures that
  $\mathfrak{C}_G=\mathfrak{Q}$.  Since $(\mathfrak{Q}, \subseteq)$ is a
  join-semilattice, Proposition\ \ref{prop:joinSL-hasse} implies that
  $\Hasse(\mathfrak{Q}, \subseteq)$ has the global
  pairwise-lca-property.  Since $G$ and $\Hasse(\mathfrak{Q},\subseteq)$
  are isomorphic, $G$ has the global pairwise-lca-property.  By
  Lemma~\ref{lem:globallca=>network}, $G$ is a network and $\mathfrak{C}_G
  = \mathfrak{Q}$ is a clustering system.
  \end{proof}
  
\begin{theorem}\label{thm:NgLCA-join.stuff}
  A network $N$ has the global lca-property if and only if
  $(V(N),\preceq_N)$ is a join-semilattice.
\end{theorem}
\begin{proof}
  One easily verifies that $\sf(N)= \hasse((V(N),\preceq_N))$.  By
  Lemma~\ref{lem:shortcutfree}, $N$ has the global lca-property if and
  only if $\sf(N) = \hasse((V(N),\preceq_N))$ has the global
  lca-property.  By Proposition~\ref{prop:joinSL-hasse}, the latter is
  equivalent to $(V(N),\preceq_N)$ being a join-semilattice.
\end{proof}

\section{A characterization by $\Ktt$-minors}
\label{sec:Ktt}

\begin{figure}
  \centering
  \includegraphics[width=0.8\textwidth]{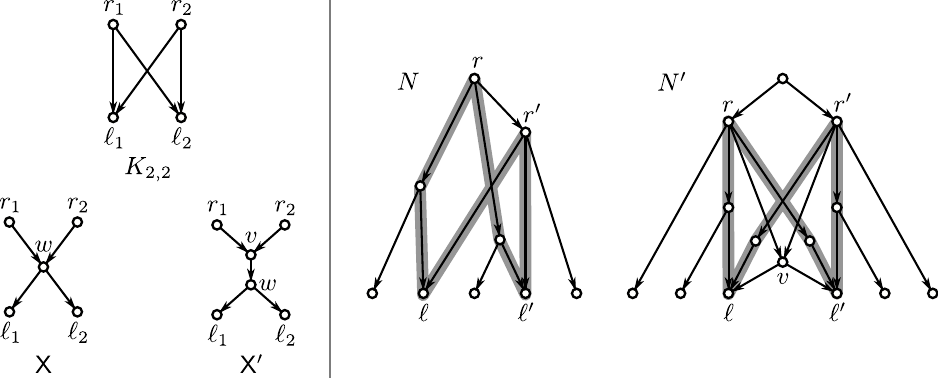}
  \caption{
  	The three DAGs $\Ktt$, {$\mathsf{X}$} and {$\mathsf{X}'$} that appear in the characterization of global lca-networks in
  	terms of minors (Thm.~\ref{thm:globalca<=>findXminor-in-K22minors}). In addition, 
  	two global lca-networks $N$ and $N'$ are shown. 
  	Here, $N$  contains a $\Ktt$-minor but not a strict $\Ktt$-minor and  $N'$ contains
  	a strict $\Ktt$-minor (the respective $\Ktt$-minors are located within the gray-shaded area).
  }\label{fig:minors}
\end{figure}

An \emph{edge-subdivision} of a graph $G$ is the replacement of an edge
$(x,y)$ in $G$ by the two edges $(x,z)$ and $(z,y)$ for some new vertex
$z\notin V(G)$.  A graph $G'$ is \emph{subdivision of $G$} if $G'$ can be
obtained from $G$ by a sequence of edge-subdivisions.  A graph $H$ is a
\emph{(topological) minor} of $G$, also known as a homeomorphic subgraph,
if a subdivision of $H$ is isomorphic to a subgraph $G'$ of $G$.  In this
case, we say that $G$ contains an $H$-minor and that $G'$ is an
$H$-subdivision.  A $\Ktt$ is the DAG having two roots $r_1$ and $r_2$, two
leaves $\ell_1$ and $\ell_2$ and edges ($r_i,\ell_j)$, $1\leq i,j\leq 2$, depicted
in Figure~\ref{fig:minors}.  At a first glance, one might assume that the
appearance of a $\Ktt$-minor in a DAG $G$ indicates the existence of two
vertices $x$ and $y$ such that $|\LCA_G(\{x,y\})|\geq2$.  However, the
network $N$ of Figure~\ref{fig:minors} exemplifies that this is not the
case; $N$ has a $\Ktt$-minor yet still has the global lca-property. We need
the following additional structure on the $\Ktt$-minors to find a
connection to unique least common ancestors.
	 
\begin{definition}\label{def:Ktt-submin}
  A DAG $G$ contains a \emph{strict $\Ktt$-minor} if it contains
  $\Ktt$-subdivision $H\subseteq G$ such that $r$ and $r'$, respectively
  $\ell$ and $\ell'$, are $\preceq_G$-incomparable for $R(H)=\{r,r'\}$ and
  $L(H)=\{\ell,\ell'\}$. In this case, we call $H$ a \emph{strict
  $\Ktt$-subdivision}.
\end{definition}

\begin{lemma}\label{lem:sufficient-K22minor}
  Let $G=(V,E)$ be a DAG. If there exists $u,v\in V$ such that
  $|\LCA_G(\{u,v\})|\geq 2$, then $G$ contains a $\Ktt$-minor. In
  particular, there is a strict $\Ktt$-subdivision $H \subseteq G$ such
  that $R(H)\subseteq\LCA_G(\{u,v\})\cap\LCA_G(L(H))$.
\end{lemma}
\begin{proof}
  Let $G=(V,E)$ be a DAG and assume $u,v\in V$ are vertices such that
  $|\LCA_G(\{u,v\})|\geq 2$.  In particular, $u\neq v$ holds and we may
  pick distinct $r,r'\in \LCA_G(\{u,v\})$. Now, since $u,v\prec_G r$ and
  $u,v\prec_G r'$ four directed paths necessarily exist in $G$: $P_{ru}$,
  $P_{rv}$, $P_{r'u}$ and $P_{r'v}$, where $P_{qw}$ denotes some path from
  $q\in\{r,r'\}$ to $w\in\{u,v\}$. Since $r$ and $r'$ are both LCAs of
  $\{u,v\}$, no vertex $p\prec_G r$ is a common ancestor of both $u$ and
  $v$, and no vertex $p'\prec_G r'$ is a common ancestor of both $u$ and
  $v$. It is not difficult to see that this implies that
  \begin{equation}\label{eq:xypaths}
    V(P_{ru})\cap V(P_{rv})=\{r\},\quad V(P_{ru})\cap
    V(P_{r'v})=\emptyset,\quad V(P_{r'u})\cap V(P_{r'v})=\{r'\}\text{, and}
    \quad V(P_{r'u})\cap V(P_{rv})=\emptyset.
  \end{equation}
  However, it is possible that the paths $P_{ru}$ and $P_{r'u}$, or
  $P_{rv}$ and $P_{r'v}$, resp., have multiple vertices in common.  Let
  $u'$ be the $\preceq_G$-maximal vertex in $V(P_{ru})\cap V(P_{r'u})$ and
  let $P_{ru'}$ and $P_{r'u'}$ be the subpaths of $P_{ru}$ respectively
  $P_{r'u}$ from $r$ respectively $r'$ to $u'$.  Similarly, let $v'$ be the
  $\preceq_G$-maximal vertex in $V(P_{rv})\cap V(P_{r'v})$ and define
  $P_{rv'}$ and $P_{r'v'}$ analogously as the respective subpaths from $r$
  to $v'$ and $r'$ to $v'$. By definition, $V(P_{ru'})\cap
  V(P_{r'u'})=\{u'\}$ and $V(P_{rv'})\cap V(P_{r'v'})=\{v'\}$ holds. This
  argument together with Equ.~\eqref{eq:xypaths} shows that the only common
  vertices of $P_{ru'}$, $P_{r'u'}$, $P_{rv'}$ and $P_{r'v'}$ are $r,r',u'$
  and $v'$. Hence, the subgraph $H$ of $G$ that is composed of the
  respective vertices and edges of the latter four paths form a subdivision
  of a $\Ktt$, i.e., $G$ contains $\Ktt$-minor. In particular
  $R(H)=\{r,r'\}\subseteq\LCA_G(\{u,v\})$ and $L(H)=\{u',v'\}$ holds by
  construction.  By definition of LCAs, $r$ and $r'$ must be
  $\preceq_G$-incomparable.  Since both $u\preceq_G u'\prec_G r$,
  $v\preceq_G v'\prec_G r$ and $r\in\LCA_G(\{u,v\})$ holds, it follows that
  $r\in\LCA_G(\{u',v'\})$. Similarly, $r'\in\LCA_G(\{u',v'\})$, hence
  $R(H)\subseteq\LCA_G(\{u,v\})\cap\LCA_G(L(H))$. In particular,
  $|\LCA_G(\{u',v'\})|\geq|\{r,r'\}|=2$. Therefore, $u'$ and $v'$ must be
  $\preceq_G$-incomparable.  In summary, $H$ is a subdivision of a $\Ktt$
  such that the vertices in $R(H)$, resp., $L(H)$ are
  $\preceq_G$-incomparable. Hence, $G$ contains a strict $\Ktt$-minor.
\end{proof}

Note that the contrapositive of the first statement of
Lemma~\ref{lem:sufficient-K22minor} implies that if $N=(V,E)$ is a network
that does not have a $\Ktt$ as a minor, then $|\LCA_N(\{u,v\})|\leq 1$ for
all $\{u,v\}\in\binom{V}{2}$. Since $|\LCA_N(\{u,v\})|>0$ holds in all
networks, we therefore have $|\LCA_N(\{u,v\})|=1$ for all
$\{u,v\}\in\binom{V}{2}$ i.e $N$ has the global
pairwise-lca-property. By Theorem~\ref{thm:global-k2}, this implies that
$\lca_N(Y)$ is well-defined for all non-empty $Y\subseteq V$. That is, we
obtain
\begin{corollary}\label{cor:noK22=>global-lca}
  If a network $N$ does not have a $\Ktt$-minor, then $N$ is a global
  lca-network.
\end{corollary}

Neither Lemma~\ref{lem:sufficient-K22minor} nor
Corollary~\ref{cor:noK22=>global-lca} can be strengthened to equivalences.
To see this, consider the global lca-network $N'$ in
Figure~\ref{fig:minors}, which contains a strict $\Ktt$-subdivision $H$
whose vertices and edges are indicated within gray-shaded area. In
particular, the leaves $\ell$ and $\ell'$ of $H$ satisfy
$\lca_{N'}(\ell,\ell')=v$ and, informally speaking, the existence of this
vertex $v$ is displayed as an {$\mathsf{X}$}-shaped pattern ``inside''
the $\Ktt$-minor.  To make this statement more precise, we let
{$\mathsf{X}$} denote the DAG with five vertices
$r_1,r_2,w,\ell_1,\ell_2$ and directed edges $(r_i,w)$ and $(w,\ell_i)$,
$1\leq i \leq 2$. Moreover, let {$\mathsf{X}'$} denote the DAG with six
vertices $r_1,r_2,v,w,\ell_1,\ell_2$ and directed edges $(v,w)$, $(r_i,v)$
and $(w,\ell_i)$, $1\leq i \leq 2$. See Figure~\ref{fig:minors} for
depictions of the DAGs ${\mathsf{X}}$ and ${\mathsf{X'}}$.

\begin{theorem}\label{thm:globalca<=>findXminor-in-K22minors}
  A network $N$ has the global lca-property if and only if, for every
  strict $\Ktt$-subdivision $H\subseteq N$, there exists an
  {$\mathsf{X}$}- or {$\mathsf{X'}$}-subdivision $H' \subseteq N$
  with $L(H)=L(H')$ and $R(H)=R(H')$.
\end{theorem}
\begin{proof}
  Let $N=(V,E)$ be a network. First assume that $N$ has the global
  lca-property. If $N$ does not contain any strict $\Ktt$-minor the
  statement is vacuously true. Hence assume that $N$ contains a strict
  $\Ktt$-minor. Let $H\subseteq N$ be a strict $\Ktt$-subdivision and let
  $R(H)=\{r,r'\}$ and $L(H)=\{\ell,\ell'\}$.  By definition of strict
  $\Ktt$-minors, $r$ and $r'$, resp., $\ell$ and $\ell'$ are
  $\preceq_N$-incomparable.  Since $N$ has the global lca-property,
  $w\coloneqq \lca_N(\ell,\ell')$ is a well-defined vertex satisfying
  $\ell,\ell'\prec_N w$, where $w\neq \ell,\ell'$ follows from the fact
  that $\ell$ and $\ell'$ are $\preceq_N$-incomparable.  Hence two directed
  paths $P_{w\ell}$ and $P_{w\ell'}$ exist in $N$, from $w$ to $\ell$ and
  $\ell'$, respectively. Since $w$ is the least common ancestor of $\ell$
  and $\ell'$, it must hold that $V(P_{w\ell})\cap V(P_{w\ell'})=\{w\}$.
  Since both $r$ and $r'$ are $\preceq_N$-incomparable common ancestors of
  $\{\ell, \ell'\}$, Lemma~\ref{lem:common-anc-above-lca} implies that
  $w\prec_N r$ and $w\prec_N r'$.  Therefore, two directed paths $P_{rw}$
  and $P_{r'w}$ from $r$ respectively $r'$ to $w$ exist in $N$. Let $v$
  denote the $\preceq_N$-maximal vertex contained in $V(P_{rw})\cap
  V(P_{r'w})$. Note that $w\preceq_N v\prec_N r,r'$. Moreover, let
  $P_{r'v}$ denote the directed subpath of $P_{r'w}$ from $r'$ to $v$. By
  construction, the four directed paths $P_{rw}$, $P_{r'v}$, $P_{w\ell}$
  and $P_{w\ell'}$ are edge-disjoint. It is now a trivial exercise to show
  that the subgraph of $N$ consisting of the vertices and edges of these
  four paths is a subdivision of {$\mathsf{X}$} whenever $w=v$ and a
  subdivision of {$\mathsf{X'}$} whenever $w\prec_N v$. In other words,
  $N$ contains a {$\mathsf{X}$}- or {$\mathsf{X'}$}-subdivision
  $H'$ which satisfy $L(H')=\{\ell, \ell'\}=L(H)$ and
  $R(H')=\{r,r'\}=R(H)$.

  Assume now that $N$ does not have the global lca-property. By
  Theorem~\ref{thm:global-k2}, $N$ does not have the global
  pairwise-lca-property either, i.e., there are distinct vertices $u,v\in
  V$ such that $|\LCA_N(u,v)|\neq1$. Since $N$ is a network,
  $\LCA_N(\{u,v\})$ is non-empty and we may therefore apply
  Lemma~\ref{lem:sufficient-K22minor} to find a strict $\Ktt$-subdivision
  $H\subseteq N$ that satisfies
  $R(H)\subseteq\LCA_N(\{u,v\})\cap\LCA_N(L(H))$. Assume, for
  contradiction, that there exists an {$\mathsf{X}$}- or
  {$\mathsf{X'}$}-subdivision $H'\subseteq N$ such that $L(H')=L(H)$
  and $R(H')=R(H)$.  One easily verifies that, in $H'$, {there exists a
    vertex $w$ such} that $\LCA_{H'}(L(H'))=\{{w}\}$.  However, since
  $R(H')=\{r,r'\}=R(H)$ and {$w\in V(H')$} it also holds that
  ${w}\preceq_N r, r'$. Since $r$ and $r'$ are $\preceq_N$-incomparable
  by assumption, even ${w}\prec_N r, r'$ holds
  (cf.\ Lemma~\ref{lem:common-anc-above-lca}).  In particular, ${w}$ is
  a common ancestor of $L(H)=L(H')$.  Moreover, since $r,r'\in
  \LCA_N(L(H))$ {the latter arguments show that $w$ is a common
    ancestor of $L(H)$} strictly below vertices contained in
  $\LCA_N(L(H))$; a contradiction.
\end{proof}

\section{A constructive characterization}
\label{sec:holju}

In this section, we provide a way to construct networks with
global lca-property starting from the graph $K_1$ that consists of a single
vertex and no edges. To this end, we show first that every network can be
constructed by stepwise adding leaves. We will later refine this process
to obtain global lca-networks. To establish this result, we note that a
directed graph $G$ is a DAG if and only if it admits a topological order
$\lessdot_G$ defined on $V(G)$, i.e., the vertices of $G$ can be totally
ordered such that $(x,y)\in E(G)$ implies that $x\lessdot_G y$
\cite{Kahn:62}.

\begin{lemma}\label{lem:network-construction}
  Let $\mathbf{N}$ denote the class of all networks. Then, $N\in
  \mathbf{N}$ if and only if $N\simeq K_1$ or $N$ can be constructed from a
  network $N'\in \mathbf N$ by applying the following operation:
  \begin{enumerate}[label=(O)]
  \item Add a leaf $x\notin V(N')$ and edges $(w,x)$ for all $w\in W$, for
    some non-empty {subset} $W\subseteq V(N)$. \label{rule:add-x-N}
  \end{enumerate}
\end{lemma}
\begin{proof}
  Let $N\in \mathbf{N}$. If $N\simeq K_1$, then we are done. Hence, assume
  that $N$ contains at least two vertices. It is easy to verify that
  $L(N)\neq \emptyset$ and $\rho_N\notin L(N)$
  (cf.\ e.g.\ \cite[p.\ 6]{HL:24}). Let $x\in L(N)$ and put $W =
  \parent_N(x)$. Consider now the subgraph $N'$ of $N$ obtained from $N$ by
  removal of $x$ and its incident edges. As $x\neq\rho_N$ is a leaf and the
  only vertex that has been removed, $N'$ remains a DAG with unique root
  $\rho_N$, i.e., $N'\in \mathbf{N}$. Applying Rule~\ref{rule:add-x-N} to
  $N'$ w.r.t.\ the set $W$ yields now $N$, which completes the proof of the
  \emph{only-if} direction.
	
  For the \emph{if} direction assume that $N\simeq K_1$ or $N$ can be
  constructed from a network $N'\in \mathbf N$ by applying
  Rule~\ref{rule:add-x-N}. Again if $N\simeq K_1$, then $N\in
  \mathbf{N}$. Assume that $N\not\simeq K_1$. Since $N'$ is, in particular,
  a DAG, there is a topological order $\lessdot_{N'}$ defined on
  $V(N')$. We can now extend this $\lessdot_{N'}$ to a topological order
  $\lessdot_{N}$ defined on $V(N)$ by simply adding $x$ as the
  $\lessdot_{N}$-maximal element. Thus, $N$ is a DAG
  \cite{Kahn:62}. Moreover, $N$ still has only one single root, since the
  construction of $N$ ensures that $\indeg_N(v)=0$ if and only if $v$
  coincides with the unique root of $N'$. The latter two arguments imply
  that $N\in \mathbf{N}$.
\end{proof}

Suppose that $N$ is obtained from $N'$ by applying rule \ref{rule:add-x-N},
i.e., we add a new leaf $x$ to $N'$ and connect $x$ to the vertices in
$\emptyset\neq W\subseteq V(N')$, then $N' =N -x \coloneqq (V(N)\setminus
\{x\}, E(N)\setminus\{(w,x)\mid w\in W\})$.  Hence, it is easy to verify
that the following property holds
\begin{enumerate}[label=(P\arabic*)]  
\item $u\preceq_{N'}v$ if and only if $u\preceq_{N} v$, for all $u,v\in
  V(N')\subseteq V(N)=V(N')\cup\{x\}$. \label{property:prec}
\end{enumerate}
In particular, \ref{property:prec} implies that the following property is
satisfied:
\begin{enumerate}[label=(P\arabic*)]  \setcounter{enumi}{1}
\item $\LCA_N(A)=\LCA_{N'}(A)$ for all subsets $A\subseteq
  V(N')$. \label{property:lca}
\end{enumerate}
Hence, if $N' = N-x$ does not have the global lca-property, then there is
some subset $A\subseteq V(N')$ such that $|\LCA_{N'}(A)| \neq 1$. This
together with \ref{property:lca} implies the following
\begin{observation}\label{obs:glblca}
  If $x$ is a leaf in a network $N$ and $N' = N-x$ does not have the
  global lca-property, then $N$ does not have the global lca-property.
\end{observation}

\begin{figure}
  \centering
  \includegraphics[width=0.8\textwidth]{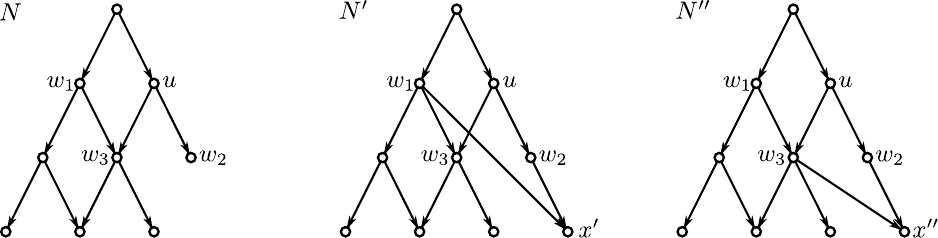}
  \caption{Three networks $N$, $N'$ and $N''$, where $N'$ is obtained from $N$ by applying
           \ref{rule:add-x-N} but not \ref{rule:addLeaf+}, and $N''$ is obtained from
           $N$ by applying \ref{rule:addLeaf+}. Both $N$ and $N''$ have the global lca-property,
           while $N'$ does not. See the text for more details.}
  \label{fig:O*-operation}
\end{figure}

By Lemma~\ref{lem:network-construction}, every network that has the
global lca-property can be constructed starting from a $K_1$ and repeated
applications of Rule~\ref{rule:add-x-N}. However, not all networks are
global lca-networks. Hence, to provide a construction of global 
lca-networks, we must refine Rule~\ref{rule:add-x-N} and, in particular,
restrict the type of sets $W$ used to add edges $(w, x)$ to {$N'$}. To
obtain this result, the following set
\begin{equation}\label{eq:L(Wv)-LCA}
  \mathcal{L}_{N'}(W|v) \coloneqq \bigcup_{w\in W}\LCA_{{N'}}(\{w,v\})
\end{equation}
defined for networks $N'$, a subset $W\subseteq V(N')$ and a given vertex
$v\in V(N')$, will play a central role.  In essence,
$\mathcal{L}_{N'}(W|v)$ is the set comprising all LCAs of $v$ and $w$
for all $w\in W$ in $N'$. For example, considering the network $N$ in
Figure~\ref{fig:O*-operation}, we have
$\mathcal{L}_N(\{w_1,w_2\}|w_3)=\LCA_N(\{w_1,w_3\})\cup
\LCA_N(\{w_2,w_3\})=\{w_1,u\}$.  Note that if $N'$ has
global lca-property, then $|\LCA_{{N'}}(\{w,v\})| = 1$ for all $v,w\in V(N')$
and Equation~\eqref{eq:L(Wv)-LCA} simplifies to
\begin{equation}\label{eq:L(Wv)-lca}
\mathcal{L}_{N'}(W|v) \coloneqq \{\lca_{N'}(v,w)\,:\, w\in W\}. 
\end{equation}

\begin{lemma}\label{lem:blopp}
  Let $N'$ be a graph with the global lca-property, and assume that $N$
  is obtained from $N'$ by applying \ref{rule:add-x-N}, where $W =
  \parent_N(x)$ denotes the set of parents for the newly added leaf
  $x\notin V(N')$. Then $N$ has the global lca-property if and only if,
  for all $v\in V(N')$, the set $\mathcal{L}_{N'}(W|v)$ contains a unique
  $\preceq_{N'}$-minimal vertex.
\end{lemma}
\begin{proof}
  Let $N'$ be a network with the global lca-property and suppose that
  $N$ is obtained from $N'$ by application of \ref{rule:add-x-N} for some
  non-empty $W\subseteq V(N')$, attaching the new leaf $x$ with the edges
  $(w,x)$, for each $w\in W$. By construction, we have $W =
  \parent_N(x)$. Note that, by \ref{property:lca}, $\LCA_N(A)=\LCA_{N'}(A)$
  for all subsets $A\subseteq V(N')$. Since $N'$ has the
  global lca-property, it follows that, in particular,
  \begin{enumerate}[label=(P\arabic*$^*$)]  \setcounter{enumi}{1}
  \item $\lca_N(u,v)=\lca_{N'}(u,v)$ is well-defined for all $u,v\in
    V(N')$. \label{property:lca-extra}
  \end{enumerate}
  
  Before showing the equivalence of the statement of this lemma, we
  consider the set $\mathcal{L}^{\min}_{N'}(W|v)$ that comprises all
  $\preceq_{N'}$-minimal vertices in $\mathcal{L}_{N'}(W|v)$ and show that,
  for every vertex $v\in V(N')$, it holds that
  \begin{equation}\label{eq:LCAsetW}
    \LCA_N(x,v) = \mathcal{L}^{\min}_{N'}(W|v), 
  \end{equation}
  that is, the set of LCAs of $v$ and $x$ in $N$ coincide with the set
  $\mathcal{L}^{\min}_{N'}(W|v)$.  Note that $\mathcal{L}_{N'}(W|v)$ can be
  chosen as specified in Equation~\eqref{eq:L(Wv)-lca}.
  
  Let $v\in V(N')=V(N)\setminus\{x\}$ be chosen arbitrarily. We show first
  that $\LCA_N(x,v) \subseteq \mathcal{L}^{\min}_{N'}(W|v)$. Note, by
  Lemma~\ref{lem:network-construction}, $N$ is a network and thus,
  $\LCA_N(x,v)\neq \emptyset$. Hence, we can arbitrarily choose a vertex
  $u\in \LCA_N(x,v)$. Since $x$ is a leaf of $N$ and $v\neq x$, we have
  $x\prec_N u$ and there exists some $p\in\parent_N(x)= W$ such that
  $x\prec_N p \preceq_N u$. In particular, $u\in V(N')$ and
  \ref{property:prec} ensures that $v\preceq_{N'}u$ and $p\preceq_{N'}u$
  i.e.\ $u$ is a common ancestor of $v$ and $p$ in $N'$. Since, $N'$ has
  the global lca-property, $\lca_{N'}(v,p)$ is well-defined and, by the
  latter arguments, satisfies $\lca_{N'}(v,p)\preceq_{N'}u$. Since
  $u,v,p\in V(N')$, we can apply \ref{property:prec} and
  \ref{property:lca-extra} to conclude that $z\coloneqq
  \lca_{N'}(v,p)=\lca_N(v,p) \preceq_{N}u$. Since $x\prec_N p\preceq_Nz$
  and $v\preceq_N z$, the vertex $z$ is a common ancestor of $v$ and $x$ in
  $N$. Together with $u\in \LCA_N(x,v)$ {the latter} implies that $z
  \not\prec_{N}u$.  Consequently, $u = z =\lca_{N'}(v,p)$ and, therefore,
  $u\in\mathcal{L}_{N'}(W|v)$. It remains to show that $u\in
  \mathcal{L}^{\min}_{N'}(W|v)$. Assume, for contradiction, that $u$ is not
  a $\preceq_{N'}$-minimal vertex in $\mathcal{L}_{N'}(W|v)$. Hence, there
  exists some vertex $w\in W$ such that $\lca_{N'}(v,w)\prec_{N'} u$. Since
  $u,v,w \in V(N')$, we can use \ref{property:prec} and
  \ref{property:lca-extra} to conclude that $\lca_{N'}(v,w) = \lca_{N}(v,w)
  \prec_{N}u$. Since $w\in W = \parent_N(x)$, $\lca_{N}(v,w)$ is also a
  common ancestor of $x$ and $v$ in $N$. This and $\lca_{N}(v,w)
  \prec_{N}u$ contradicts the fact that $u\in\LCA_N(x,v)$. Thus, $u\in
  \mathcal{L}^{\min}_{N'}(W|v)$. As $v\in V(N')$ and $u\in
  {\LCA_N(x,v)}$ have been chosen arbitrarily it follows that
  $\LCA_N(x,v) \subseteq \mathcal{L}^{\min}_{N'}(W|v)$ {for all $v\in
    V(N')$}.

  We proceed with showing that $\mathcal{L}^{\min}_{N'}(W|v) \subseteq
  \LCA_N(x,v)$. Since $N'$ has the global lca-property and since $W\neq
  \emptyset$, it follows that $\mathcal{L}^{\min}_{N'}(W|v)\neq
  \emptyset$. Let $u\in \mathcal{L}^{\min}_{N'}(W|v)$ be chosen arbitrarily
  and let $w\in W$ be a vertex such that $\lca_{N'}(v,w)=u$. Since $v,w\in
  V(N')$, \ref{property:lca-extra} implies that $u = \lca_{N}(v,w)$. Since
  $w\in W = \parent_N(x)$, $u$ is a common ancestor of $x$ and $v$ in
  $N$. Hence, there is a vertex $u'\in\LCA_N(x,v)$ such that $u'\preceq_N
  u$. To recall, $x$ is a leaf of $N$ and $v\neq x$ which implies that
  $x\prec_N u'$ must hold and there exists some parent $p\in W$ of $x$ in
  $N$ such that $x\prec_N p\preceq_N u'$.  This and $u'\in\LCA_N(x,v)$
  implies that $u'$ is a common ancestor of $v$ and $p$ in $N$. By
  \ref{property:prec}, $u'$ is a common ancestor of $v$ and $p$ in
  $N'$. Since $N'$ has the global lca-property, $\lca_{N'}(v,p)$ is
  well-defined. Taking the latter two arguments together we obtain
  $\lca_{N'}(v,p)\preceq_{N'} u'\preceq_{N'}u$. By construction,
  $\lca_{N'}(v,p)\in \mathcal{L}_{N'}(W|v)$. Moreover, since $u\in
  \mathcal{L}^{\min}_{N'}(W|v)$, the vertex $u$ is a $\preceq_{N'}$-minimal
  vertex in $\mathcal{L}_{N'}(W|v)$ and thus, $\lca_{N'}(v,p) \prec_{N'} u$
  cannot hold. Taking the latter arguments together we have $\lca_{N'}(v,p)
  = u' = u$ which, in turn, together with $u'\in\LCA_N(x,v)$ implies that
  $u\in\LCA_N(x,v)$. As $u\in \mathcal{L}^{\min}_{N'}(W|v)$ was chosen
  arbitrarily, it follows that $ \mathcal{L}^{\min}_{N'}(W|v) \subseteq
  \LCA_N(x,v)$. In summary, we have shown that
  $\mathcal{L}^{\min}_{N'}(W|v) = \LCA_N(x,v)$ {for all $v\in V(N')$}.

  We now return to the main statement. Suppose first that $N$ has the
  global lca-property. In this case, Theorem~\ref{thm:global-k2} implies
  that $N$ has the global pairwise-lca-property. This is, in particular,
  only possible if $|\LCA_N(x,v)|=1$ for all $v\in V(N')\subseteq V(N)$. As
  $ \mathcal{L}^{\min}_{N'}(W|v) = \LCA_N(x,v)$ for all $v\in V(N')$ it
  follows that $|\mathcal{L}^{\min}_{N'}(W|v)|=1 $ and, therefore, that
  $\mathcal{L}_{N'}(W|v)$ contains a unique $\preceq_{N'}$-minimal element
  for all $v\in V(N')$. Conversely, suppose that $\mathcal{L}_{N'}(W|v)$
  contains a unique $\preceq_{N'}$-minimal element for all $v\in V(N')$ and
  thus that $|\mathcal{L}^{\min}_{N'}(W|v)|=1 $ for all $v\in V(N')$. In
  this case, $ \mathcal{L}^{\min}_{N'}(W|v) = \LCA_N(x,v)$ implies that
  $|\LCA_N(x,v)|=1$ for all $v\in V(N')$.  Clearly, $|\LCA_N(x,x)|=1$. In
  addition, $|\LCA_N(a,b)|=1$ for all $a,b\in V(N')$ as $N'$ has the
  global lca-property (c.f.\ \ref{property:lca-extra}). The latter three
  arguments imply that $|\LCA_N(a,b)|=1$ for all $a,b\in V(N)$. Thus, $N$
  has the global pairwise-lca-property and Theorem~\ref{thm:global-k2},
  implies that $N$ has the global lca-property, which completes this
  proof.
\end{proof}

In the following, we show that the class of networks with the
global lca-property coincides with the class of \holju
graphs\footnote{We call these graphs \holju{} because the main concept was
conceived by some of us while being stranded at Ljubljana Airport (LJU)
with a rather significant hangover (HO). We thank the organizers of the
40th TBI-winterseminar in Bled, Slovenia for providing an endless amount of
free beer, schnaps and wine, and assure the reader that all proofs in this
manuscript were written after a full dose of sleep.}
\begin{definition}\label{def:holju}
  The class of \holju graphs is defined recursively.  The graph $K_1$ is a
  \holju graph. If $G$ is a \holju graph, then any graph obtained by
  applying the following operation once, is also a \holju graph:
  \begin{enumerate}[label= (O$^*$)] 
  \item Add a leaf $x\notin V(G)$ and the edges $(w,x)$ for all $w\in W$,
    where $W\subseteq V(G)$ is a non-empty subset of vertices such that for
    any $v\in V(G)$ the set $\mathcal{L}_{G}(W|v)$ contains a unique
    $\preceq_G$-minimal vertex. Note that the case $|W|=1$ is
    allowed. \label{rule:addLeaf+}
    \end{enumerate}
\end{definition}
Note that \ref{rule:addLeaf+} is just a special case of \ref{rule:add-x-N}. For example,
consider the three networks $N$, $N'$ and $N''$ in Figure~\ref{fig:O*-operation}, where
$N=N'-x'=N''-x''$. Note that $N$ has the global lca-property, and that both $N'$ and $N''$ is
obtained from $N$ by applying \ref{rule:add-x-N}, although to different sets $W$. For $N'$, we have
$W=\{w_1,w_2\}$ and, in particular, $\mathcal{L}_N(\{w_1,w_2\}|w_3)=\{w_1,u\}$. Since $w_1$ and $u$ are
$\preceq_N$-incomparable, $\mathcal{L}_N(\{w_1,w_2\}|w_3)$ does not contain a unique $\preceq_N$-minimal vertex
and, consequently, $N'$ is not obtained from $N$ by applying \ref{rule:addLeaf+}. In contrast, $N''$
is obtained from $N$ by applying \ref{rule:addLeaf+} to $W=\{w_2,w_3\}$, realized by (tediously but
with no great effort) verifying that $\mathcal{L}_N(W|x)$ indeed has a unique $\preceq_N$-minimal
vertex for every $x\in V(N)$. In this example, $N'$ does not have the global lca-property (since
e.g. $|\LCA_{N'}(\{x',w_3\})|=|\{w_1,u\}|\neq 1$) while $N''$ is a global lca-network. This is not a
coincidence, due to the following theorem.

\begin{theorem}\label{thm:hlojy}
  A network has the global lca-property if and only if it is a \holju
  graph.
\end{theorem}
\begin{proof}
  Let $N$ be a \holju graph. Since \ref{rule:addLeaf+} is just a special
  case of \ref{rule:add-x-N}, Lemma~\ref{lem:network-construction} ensures
  $N$ is a network. Note that the constraints imposed on
  $\mathcal{L}_{G}(W|v)$ in Lemma~\ref{lem:blopp} and \ref{rule:addLeaf+}
  are the same. As $N$ is constructed by a sequence of operations of
  \ref{rule:addLeaf+} starting with the $K_1$ that trivially has the
  global lca-property, repeated application of Lemma \ref{lem:blopp},
  therefore, implies that $N$ has the global lca-property.
  
  Suppose now that $N$ is a network with the global lca-property. If
  $N\simeq K_1$, then $N$ is clearly a \holju graph. Suppose that
  $N\not\simeq K_1$. Since $N$ is a network, it can, due to
  Lemma~\ref{lem:network-construction}, be constructed by stepwise
  application of \ref{rule:add-x-N} starting with a $K_1$. This yields a
  sequence of constructed networks $K_1 = N_1,\, N_2,\, \dots,\, N_k = N$,
  where $k>1$. Let us denote with $x_i$ the new leaf added to $N_{i-1}$ and
  with $W_{i}$ the set of parents of $x_i$ in $N_i$ used in the
  construction of $N_{i}$ when applying \ref{rule:add-x-N} to $N_{i-1}$,
  for each $i\in \{2,\dots,k\}$. In other words, $N_{i-1} = N_i - x_i$ for
  the leaf $x_i$ of $N_i$, for each $i\in \{2,\dots,k\}$. Since $N=N_k$ has
  the global lca-property, Observation~\ref{obs:glblca} implies that
  $N_{k-1}$ must have the global lca-property. By induction, $N_i$ must
  have the global lca-property, $1\leq i\leq k$. Since each $N_i$ is
  obtained from $N_{i-1}$ by applying \ref{rule:add-x-N} and since
  $N_{i-1}$ has the global lca-property, Lemma \ref{lem:blopp}
  {implies that}
  $\mathcal{L}_{N_{i-1}}(W_i|v)$ contains a unique
  $\preceq_{N_{i-1}}$-minimal vertex, for each $i\in \{2,\dots,k\}$. Hence,
  to obtain $N_i$ from $N_{i-1}$ we have, in fact, applied
  \ref{rule:addLeaf+}, $i\in \{2,\dots,k\}$. Therefore, each $N_i$ and
  thus, $N$ is a \holju graph.
  \end{proof}

\section{A characterization by tree-leaf-child and  leaf-extended DAGs}
\label{sec:lxtG}

A \emph{tree-leaf} $v$ in a DAG $G$ is a leaf $v\in L(G)$ with in-degree
one.  A DAG is \emph{tree-leaf-child} if every inner vertex is adjacent to
a tree-leaf.  Tree-leaf-child DAGs form a subclass of the \emph{tree-child}
DAGs, in which every inner vertex is adjacent to vertex with with in-degree
one, but not necessarily to a leaf \cite{CRV:09}.  For later reference we
provide
\begin{lemma}[{\cite[Cor.~20]{Hellmuth:22q}}]\label{lem:TR-N-ClosedCLS}
  A tree-child network $G$ is an lca-network if and only if its clustering
  system $\mathfrak{C}_N$ is closed.
\end{lemma}

\begin{lemma}\label{lem:ltc-cluster-PCC}
  Every  tree-leaf-child DAG $G$ satisfies (PCC).
\end{lemma}
\begin{proof}
  To see that $G$ satisfies (PCC), let $u,v\in V(G)$ be two distinct
  vertices. If $u\preceq_G v$, then $\CC_{G}(u)\subseteq \CC_{G}(v)$ by
  Lemma~\ref{lem:prec-subset}. Assume now that $\CC_{G}(u)\subseteq
  \CC_{G}(v)$. If $u$ is a leaf in $G$, then $\CC_{G}(u)=\{u\}$ and thus,
  $u\in \CC_{G}(v)$. By definition, $u\preceq_{G} v$. Assume now that $u$
  is an inner vertex in $G$. Since $G$ is tree-leaf-child, there is a leaf
  $x_u$ whose unique parent is $u$. Moreover, $x_u\in \CC_{G}(v)$. Thus,
  there is a directed path $P$ from $v$ to $x_u$. Taking the latter
  arguments together, this path $P$ must contain $u$. Hence, $u\preceq_{G}
  v$. In summary, $G$ satisfies (PCC).
\end{proof}

\begin{figure}
  \centering
  \includegraphics[width=0.8\textwidth]{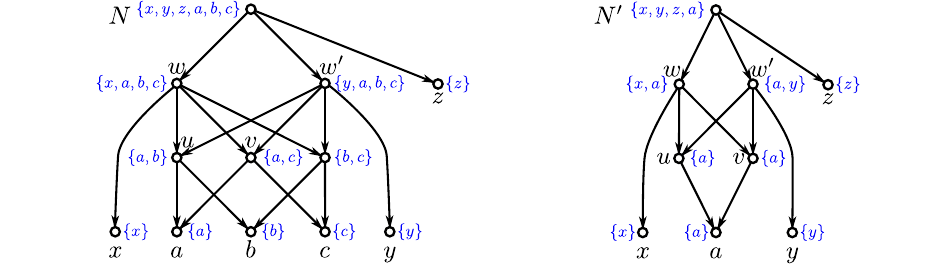}
  \caption{ Two networks $N$ and $N'$, neither of which have the global
    (pairwise)-lca-property. Moreover, $N$ does have the
    pairwise-lca-property i.e.\ $\lca_N(A)$ is well-defined for each subset
    $A\subseteq L(N)$ of size $|A|=2$, yet $\LCA_N(\{a,b,c\})=\{w,w'\}$ and
    $N$ therefore does not have the lca-property. In $N'$, there is a
    unique LCA for every non-empty subset of leaves, but
    $\LCA_{N'}(\{u,v\})=\{w,w'\}$, implying that $N'$ is not a global
    lca-network. Moreover, both $N$ and $N'$ are networks in which every
    inner vertex is adjacent to a leaf, indicating the need to enforce
    \emph{leaf-tree-child} in Proposition~\ref{prop:TLC-lca-char}. Lastly,
    for $N$, the set system $\mathfrak{C}_N$ is a pre-binary clustering
    system that is not closed since $\CC_N(w)\cap \CC_N(w') = \{a,b,c\}
    \notin \mathfrak{C}_N$. The latter argument shows that the Statements
    (5) and (6) in Proposition~\ref{prop:TLC-lca-char} are, in the general
    setting, not equivalent.}
  \label{fig:global-subtle}
\end{figure}

In general DAGs, the pairwise-lca-property and the lca-property do not
imply knowledge about the global lca-property nor of the global
pairwise-lca-property. To see this consider the network $N$ and $N'$ in
Figure~\ref{fig:global-subtle}. Here, $\lca_N(A)$ is well-defined for each
subset $A\subseteq L(N)$ of size $|A|=2$ and $N$ consequently has the
pairwise-lca-property. Despite this, it does not have the lca-property (and
therefore not the \emph{global} lca-property) since, for example,
$|\LCA_N(\{a,b,c\})|=|\{w,w'\}|\neq1$. In the network $N'$ of the same
figure, $\lca_{N'}(A)$ is well-defined for each non-empty $A\subseteq
L(N')$ and therefore $N'$ has both the lca-property and the
pairwise-lca-property, yet $|\LCA_{N'}(u,v)|=|\{w,w'\}|\neq1$ and $N'$ does
not have the global lca-property. Nevertheless, we can connect the global
(pairwise-)lca-property with the (pairwise-)lca-property in the case of
tree-leaf-child DAGs. To this end we will need to consider the following
property of clustering systems:
\begin{definition} \cite{Barthelemy:08}
  A clustering system $\mathfrak{C}$ is \emph{pre-binary} if for every
  $x,y\in X$ there is unique inclusion-minimal element $C\in\mathfrak{C}$
  that contains $x$ and $y$.
\end{definition}
We are now in the position to extend Lemma~\ref{lem:TR-N-ClosedCLS} in the
following manner:
\begin{proposition}\label{prop:TLC-lca-char}
  The following statements are equivalent for every tree-leaf-child DAG $G$:
  \begin{enumerate}[noitemsep, label=(\arabic*)]
  \item $G$ has the global lca-property.  
  \item $G$ has the global pairwise-lca-property.  
  \item $G$ has the lca-property.
  \item $G$ has the pairwise-lca-property.
  \item $\mathfrak{C}_{G}$ is a closed clustering system.
  \item $\mathfrak{C}_{G}$ is a pre-binary clustering system.
  \end{enumerate}
\end{proposition}
\begin{proof} 
  Let $G$ be a tree-leaf-child DAG on $X$. By Theorem~\ref{thm:global-k2},
  Condition (1) and (2) are equivalent.  Moreover, Condition (1) clearly
  implies (3), and Condition (3) implies (4).  We continue with showing
  that Condition (4) implies (2) to establish the equivalence between
  Conditions (1), (2), (3) and (4).

  To this end, we assume, by contraposition, that $G$ does not have the
  global pairwise-lca-property. In this case, there are distinct vertices
  $x,y\in V(G)$ such that $|\LCA_{G}(\{x,y\})|\neq 1$. If $x,y\in L(G)$,
  then it immediately follows that $G$ does not have the
  pairwise-lca-property. Hence consider the case when at least one of $x$
  and $y$ is an inner vertex of $G$. Without loss of generality, assume
  that $x\notin X$. Since $G$ is tree-leaf-child, we can find a tree-leaf
  $x'$ such that $x'$ is a child of $x$. Similarly, if $y\notin X$ we can
  find a tree-leaf $y'$ such that $y'$ is a child of $y$. If $y$ is a leaf
  itself, put $y'=y$. Since $x'\preceq_{G}x$ and $y'\preceq_{G}y$, every
  common ancestor of $x$ and $y$ is also a common ancestor of $x'$ and
  $y'$. Moreover, since $x'$ is a tree-leaf and a child of $x$, every
  ancestors of $x'$ -- except $x'$ itself -- is also an ancestor of $x$. If
  $y'\prec_{G} y$, analogous arguments shows that every ancestors of $y'$
  -- except $y'$ itself -- is also an ancestor of $y$ and, if $y=y'$, the
  ancestors of $y$ and $y'$ clearly coincide. In summary so-far, we have
  \begin{equation*}
    \ANC_G(\{x,y\})\subseteq\ANC_G(\{x',y'\})\subseteq
    \ANC_G(\{x,y\})\cup\{x',y'\}.
    \end{equation*}
  Note that $|\LCA_{G}(\{x,y\})|\neq 1$ implies that neither $x\preceq_G y$
  nor $y\preceq_G x$ can hold. Together with $x\neq y$, this implies that
  $x$ and $y$ are $\preceq_G$-incomparable.  Since $x'$ has as unique
  parent $x$ and either $y=y'$ or $y'$ has as unique parent $y$, it follows
  that $x'$ and $y'$ must be $\preceq_G$-incomparable.  It is now easy to
  see that each vertex that is a common ancestor of $x'$ and $y'$ is a
  common ancestor of $x$ and $y$, i.e., $\ANC_G(\{x,y\})
  =\ANC_G(\{x',y'\})$.  In particular
  $\LCA_{G}(\{x,y\})=\LCA_{G}(\{x',y'\})$ and, consequently, $G$ does not
  have the pairwise-lca-property.
  
  Thus, Conditions (1), (2), (3) and (4) are equivalent.  Moreover, since
  $G$ is a tree-leaf-child DAG, Lemma~\ref{lem:ltc-cluster-PCC} implies
  that $G$ satisfied (PCC). This allows us to use
  Proposition~\ref{prop:PPC+Closed=>globalLCA}, to conclude that (1) and
  (5) are equivalent.
  
  To complete the proof, we show the equivalence between Condition (4) and
  (6). First observe that Theorem~3.6 of \cite{HL:24a} states that if $H$
  is a DAG that satisfy (PCC), then $H$ has the pairwise-lca-property if
  and only if $\mathfrak{C}_{H}$ is pre-binary. As
  Lemma~\ref{lem:ltc-cluster-PCC} ensures that $G$ satisfy (PCC), the
  latter argument immediately implies that $G$ satisfies Condition (4) if
  and only if it satisfies (6).
\end{proof}

The class of tree-leaf-child DAGs is of course quite
restrictive. Nevertheless, we can associate a type of canonical
tree-leaf-child DAG $\lxt(G)$ to every DAG $G$ which will be used later to
characterize DAGs with global lca-property.
\begin{definition}
 For every DAG $G=(V,E)$, the corresponding \emph{leaf-extended DAG}
 $\lxt(G)=(V^*,E^*)$ is obtained by putting first $V^*=V$ and $E^*=E$ and
 then, adding, for every inner vertex $v\in V$, a new vertex $x_v$ to $V^*$
 and the edge $(v,x_v)$ to $E^*$.
\end{definition}
We collect now a couple of results for leaf-extended DAGs. Note that, by
definition, $\lxt(G)$ is tree-leaf-child for every DAG $G$.

\begin{lemma}\label{lem:G*-basic-properties}
  Let $G$ be a DAG and $\lxt(G)$ be its leaf-extended version. Then,
  $\outdeg_{\lxt(G)}(v)>1$ for all inner vertices in $\lxt(G)$.  Moreover,
  $\lxt(G)$ satisfies (PCC) and $\CC_{\lxt(G)}(u)\neq \CC_{\lxt(G)}(v)$ for
  all distinct $u,v\in V(\lxt(G))$.  In particular, $\lxt(G)$ is
  lca-relevant and $\sf(\lxt(G))$ is regular.
\end{lemma}
\begin{proof}
  Let $G$ be a DAG and $G^*\coloneqq \lxt(G)$ be its leaf-extended
  version. By construction, for each inner vertex $v$ in $G^*$ there is a
  unique leaf $x_v$ whose unique parent is $v$ and such that $x_v\notin
  V(G)$. In particular, $v$ is an inner vertex in $G$. Hence,
  $\outdeg_G(v)\geq 1$ and thus, $\outdeg_{G^*}(v) = \outdeg_G(v) +
  1>1$. Since $G^*$ is tree-leaf-child, Lemma~\ref{lem:ltc-cluster-PCC}
  implies that $G^*$ satisfies (PCC).
	
  To show that $G^*$ is lca-relevant, we equivalently show that the
  condition in Theorem~\ref{thm:char-lcaRel}\ref{eq:stronger-than-PCC} is
  satisfied. Suppose that $\CC_{G^*}(u)\subseteq\CC_{G^*}(v)$ for some
  vertices $u,v\in V(G^*)$. If $u$ is a leaf of $G^*$, then
  $\CC_{G^*}(u)=\{u\}$ and $u\in\CC_{G^*}(v)$ immediately implies that
  $u\preceq_{G^*}v$. If $u$ is an inner vertex of $G^*$ then, since $G^*$
  is tree-leaf-child, $u$ has a leaf-child $x_u$ such that $u$ is the
  unique parent of $x_u$ in $G^*$. Clearly, $x_u\in\CC_{G^*}(u)$ and,
  consequently, $x_u\in\CC_{G^*}(v)$. In particular any directed path from
  $v$ to $x_u$ must contain $u$, and such a path must exist by
  $x_u\in\CC_{G^*}(v)$, from which $u\preceq_{G^*}v$ follows.  The latter
  together with Lemma~\ref{lem:prec-subset} implies that
  Theorem~\ref{thm:char-lcaRel}\ref{eq:stronger-than-PCC} holds and that
  $G^*$ is lca-relevant.  Moreover,
  Theorem~\ref{thm:char-lcaRel}\ref{eq:sf-reg} implies that $\sf(G^*)$ is
  regular.  Finally, Theorem~\ref{thm:char-lcaRel}\ref{eq:PCC-neqCluster}
  ensures that $\CC_{G^*}(u)\neq \CC_{G^*}(v)$ for all distinct $u,v\in
  V(G^*)$.
\end{proof}

Note that tree-leaf-child DAGs $G$ may contain inner vertices with
out-degree one and distinct vertices $u$ and $v$ with $\CC_{G}(u) =
\CC_{G}(v)$ (e.g. a rooted edge).  Hence Lemma
\ref{lem:G*-basic-properties} does, in general, not hold for arbitrary
tree-leaf-child DAGs.

\begin{theorem}\label{thm:char-global-lca-lxt}
  The following statements are equivalent for every DAG $G$ on $X$:
  \begin{enumerate}[noitemsep, label=(\arabic*)]
  \item $G$ has the global lca-property.
  \item $\lxt(G)$ has the global lca-property.  
  \item $\lxt(G)$ has the global pairwise-lca-property.  
  \item $\lxt(G)$ has the lca-property.
  \item $\lxt(G)$ has the pairwise-lca-property.
  \item $\mathfrak{C}_{\lxt(G)}$ is a closed clustering system.
  \item $\mathfrak{C}_{\lxt(G)}$ is a pre-binary clustering system.
  \end{enumerate}
\end{theorem}
\begin{proof}
  Since $\lxt(G)$ is tree-leaf-child, Proposition~\ref{prop:TLC-lca-char}
  implies that Conditions (2)--(7) are equivalent. Therefore, we only need
  to establish the equivalence between Condition (1) and (2).  Suppose that
  the DAG $G$ has the global lca-property. By definition, $\lxt(G)$
  coincides with the DAG obtained from $G$ by using operation
  \ref{rule:add-x-N} for all subsets $W = \{w\} $ with $w\in V^0(G)$.
  Since $G$ has the global lca-property, $|\LCA_G(\{v,w\})| = 1$ for all
  $v\in V(G)$ and $w\in V^0(G)$.  Hence, for all $w\in V^0(G)$, the set
  $\mathcal{L}_{G}(\{w\}|v) = \bigcup_{w'\in \{w\}}\LCA_G(\{w',v\}) =
  \LCA_G(\{w,v\})$ as defined in Equation~\eqref{eq:L(Wv)-LCA} has a unique
  element and, thus, in particular, a unique $\preceq_G$-minimal vertex.
  Theorem~\ref{thm:hlojy} implies that $G$ is a \holju graph.  The latter
  arguments together with Definition~\ref{def:holju} imply that $\lxt(G)$
  is a \holju graph.  Again, Theorem~\ref{thm:hlojy} implies that $\lxt(G)$
  has the global lca-property.  Hence, Condition (1) implies (2).

  For the converse, assume that $G$ does not have the global lca-property,
  and let $A\subseteq V(G)$ be a non-empty subset of vertices such that
  $|\LCA_G(A)|\neq 1$. By construction, $A\subseteq V(G)\subseteq
  V(\lxt(G))$ and $\ANC_{\lxt(G)}(A)=\ANC_G(A)$. Therefore,
  $|\LCA_{\lxt(G)}(A)|\neq 1$ must hold, and $\lxt(G)$ does not have the
  global lca-property. By contraposition, Condition (2) implies (1).
\end{proof}	 

Prominent examples of phylogenetic networks in the literature are
galled-trees \cite{GSL:08} and, more general, level-1 networks
\cite{Rossello:09}.  A network $N$ is a \emph{galled-tree} if each maximal
biconnected subgraph $K$ in $N$ is either a single vertex, an edge or 
{there are two vertices $u$ and $v$ in $K$ such that}  $K$ is
composed of exactly two directed $uv$-paths that only have $u$ and $v$
in common \cite{Hellmuth:22q}. In particular, galled-trees form a subclass
of \emph{level-1 networks}, i.e., networks $N$ in which each biconnected
component contains at most one vertex $v$ with $\indeg_N(v)>1$.
\begin{corollary}\label{cor:level-1=holju}
  Every level-1 network and thus every galled-tree is a global lca-network
  and thus a \holju graph.
\end{corollary}
\begin{proof}
  Let $N=(V,E)$ be a level-1 network. It is easy to see that $\lxt(N)$
  remains a level-1 network.  As shown in Lemma~49 in \cite{Hellmuth:22q}
  every level-1 network is an lca-network.  Thus, $\lxt(N)$ is an
  lca-network. By Theorem~\ref{thm:char-global-lca-lxt}, $N$ is a global
  lca-network. By Theorem~\ref{thm:hlojy}, $N$ is a \holju graph.  Since
  galled-trees are level-1 networks, in particular galled-trees are global
  lca-network and \holju{} graphs.
\end{proof}

{The characterization in Proposition~\ref{prop:TLC-lca-char} and, via leaf-extension, in Theorem~\ref{thm:char-global-lca-lxt}, shows that for tree-leaf-child DAGs the global lca-property is equivalent to the corresponding cluster system being pre-binary, and hence equivalent to closedness. It is therefore natural to place these results into the broader context of pre-$k$-ary set systems and closure theory. The next two paragraphs briefly recall this connection and explain why the closedness conditions appearing above can be viewed as a special case of a more general set-system framework.}

In \cite{Changat:19a} a property denoted by (kKC) was considered that was
defined as follows: A set system $\mathfrak{C}$ is (kKC) if, for every
non-empty $A\subseteq V$ with $|A|\le k$ there is a unique inclusion-minimal
$C\in\mathfrak{C}$ that contains $A$. Note that, for $k=2$, this reduces to
$\mathfrak{C}$ being pre-binary. This condition can be expressed equivalently
in terms of the \emph{closure function}
\begin{equation}
  \cl(A)\coloneqq
  \bigcap_{\substack{C\in\mathfrak{C}\\ A\subseteq C}} C
  \label{eq:closure}
\end{equation}
as follows: For all non-empty $A\subseteq V$ with $|A|\le k$ it holds that
$A\subseteq \cl(A)\in\mathfrak{C}$. That is, $\cl(A)$ is the unique
inclusion-minimal $C\in\mathfrak{C}$ that contains $A$ appearing in the
definition of (kKC). Set systems with property (kKC) were called
\emph{pre-$k$-ary} in \cite{S+24}.  As a generalization, \cite{HL:24a}
considers pre-$\One$-ary set systems, where {the condition} $|A|\le k$
is replaced by $|A|\in\One$ for some set $\One$ of positive
integers. Hence, $\mathfrak{C}$ is pre-$|V|$-ary in the language of
\cite{S+24} if and only if it is pre-$\One$-ary for $\One =
\{1,\dots,|V|\}$. Theorem~3.6 in \cite{HL:24a} shows that a DAG $G$ that
satisfies (PCC) has the lca-property if and only if $\mathfrak{C}_G$ it is
pre-$|V|$-ary.

For pre-$|V|$-ary set systems, $\cl(A)$ is well-defined in
Equation~(\ref{eq:closure}) for all non-empty sets $A\subseteq V$ and
coincides with the conventional definition of the closure function. In
particular, in this case, $A\in\mathfrak{C}$ if and only if $A=\cl(A)$ and
$A\ne\emptyset$.  There is a subtle -- albeit trivial -- difference between
the definition of the closure function in the literature on convex set
systems \cite{vandeVel:93} and in the setting of clustering systems: Since
the empty set does not have an interpretation as a cluster, we define
$\emptyset\notin\mathfrak{C}$, while in convexities one usually stipulates
$\emptyset\in\mathfrak{C}$. Clustering systems therefore are considered
closed if they are closed under \emph{non-empty} intersections, and we have
$A\in\mathfrak{C}$ if and only $A=\cl(A)$ and $A\ne\emptyset$. In the
setting of convexities, we have $A\cap B\in\mathfrak{C}$ for all
$A,B\in\mathfrak{C}$ and equivalently $A\in\mathfrak{C}$ iff
$A=\cl(A)$. The definition of pre-$|V|$-ary thus coincides with
$\mathfrak{C}$ being a closed set system, i.e., an abstract convexity, up
to the exclusion of the empty set. Clearly there is a 1-1 correspondence
between set systems that include and exclude the empty set, respectively,
and definitions such as the ones of closure and closedness above are easily
adjusted to the respective setting. A brief discussion of this point can
also be found in \cite{Changat:19a}.

\section{Characterization via descendants, ancestors, and intermediaries}
\label{sec:DAB-relations}

Recall that the DAG $\lxt(G)$ is obtained from a DAG $G$ by adding, for all
inner vertices $v$ in $G$ a new tree-leaf $x_v$ whose unique parent in
$\lxt(G)$ is $v$. Hence, $x_v$ is such a new leaf in $\lxt(G)$ if and only
if $x_v\in \overline L(\lxt(G)) \coloneqq L(\lxt(G))\setminus L(G)$.  The
following technical results is used to establish further characterizations
of global lca-DAGs.  In particular, we show how {the descendant system}
$\mathfrak{D}_G$ can be obtained from $\mathfrak{C}_{\lxt(G)}$ and $L(G)$
and that $(\mathfrak{C}^*,\subseteq)$, $(\mathfrak{D}_G,\subseteq)$, and
$(V(G),\preceq_G)$ are pairwise isomorphic, where $\mathfrak{C}^*$ is
obtained from $\mathfrak{D}_G$ by removal {of} all singletons.

\begin{lemma}\label{lem:D-and-C_N*}
  For a DAG $G$, the set system $\mathfrak{D}_G$ can be obtained from
  $\mathfrak{C}_{\lxt(G)}$ by applying the following steps
  \begin{enumerate}[noitemsep]
  \item[(1)] remove the singleton $\{x_v\}$ from $\mathfrak{C}_{\lxt(G)}$ for
    all $x_v\in \overline L \coloneqq L(\lxt(G))\setminus L(G)$ and,
    afterwards,
  \item[(2)] replace all $x_v\in \overline L\cap C$ in all $C\in
    \mathfrak{C}_{\lxt(G)}$ by the respective vertex $v\in V(G)$.
  \end{enumerate}
  Moreover, if $\mathfrak{C}^*$ is the set obtained from
  $\mathfrak{C}_{\lxt(G)}$ after application of Step~(1), then there is a
  1-to-1 correspondence between the elements in $\mathfrak{C}^*$ and the
  elements in $\mathfrak{D}_G$. Furthermore, the posets
  $(\mathfrak{C}^*,\subseteq)$, $(\mathfrak{D}_G,\subseteq)$ and
  $(V(G),\preceq_G)$ are isomorphic.
\end{lemma}
\begin{proof}
  Let $G$ be a DAG and $\overline{L}\coloneqq L(\lxt(G))\setminus
  L(G)$. Since $L(G)\subseteq L(\lxt(G))$, we have $\DD{x}{G}=\{x\} \in
    \desccl{G}$ for all $x\in L(G)$. Furthermore, $x_v\notin L(G)$ for all
  $x_v\in\overline{L}$ and it follows that $\{x_v\}\notin \desccl{G}$ for
  all $x_v\in\overline{L}$. Consequently, the set of all singleton sets in
  $\mathfrak{D}_G$ can be obtained from the set of all singleton sets in
  $\mathfrak{C}_{\lxt(G)}$ by applying Step~(1).

  Assume now that Step~(1) has been applied to $\mathfrak{C}_{\lxt(G)}$
  and call the resulting set $\mathfrak{C}^*$. We proceed with showing
  that $D\coloneqq \DD{v}{G} \in \mathfrak{D}_G$ with $|D|>1$ if and only
  if $D$ can be obtained from $C\coloneqq \CC_{\lxt(G)}(v) \in
  \mathfrak{C}^*$ by applying Step (2). In other words, all non-singleton
  sets in $\desccl{G}$ can be obtained from all non-singleton sets in
  $\mathfrak{C}^*$ by applying Step~(2).
	
  Suppose that $w\in C$ and thus, that $w\in L(\lxt(G))$ and
  $w\preceq_{\lxt(G)} v$. Assume that $w\in \overline{L}$, then $w=x_z$ for
  some $z\in V(G)$. In this case, $w$ will be replaced by $z$ in $C$. This
  step is correct, as in this case, $x_z\preceq_{\lxt(G)} v$ implies that
  $z \preceq_{G} v$ and thus, $z\in D$.  Assume now that $w\notin
  \overline{L}$ and thus, that $w\in L(G)$.  Then, $w\preceq_{\lxt(G)} v$
  and thus, $w \preceq_{G} v$. As Step (2) is not applied for $w$ in this
  case, $w$ remains in $D$.  Hence, the vertices in $C$ to which Step~(2) has been
  applied yield a  subset of vertices in $D$.
	
  Suppose that $w\in D$ and thus that $w \preceq_{G} v$. Assume first that
  $w\in L(G)\subseteq L(\lxt(G))$.  Hence, $w \preceq_{G} v$ and
  $w\preceq_{\lxt(G)} v$. This and the fact that $w\notin \overline{L}$
  implies that Step (2) is not applied to $w$, i.e., $w\in C$.
  Assume now that $w\in V^0(G)$. In this case, $w\notin L(\lxt(G))$, which
  implies that $w\notin C$. However, since $w$ is an inner vertex, there is
  a tree-leaf-child $x_w\in \overline L$ with $x_w\prec_{\lxt(G)}
  w$. Moreover, $w \preceq_{G} v$ implies $w\preceq_{\lxt(G)} v$ and thus,
  $x_w\in C$. In the latter case, Step (2) is applied and $x_w$ is replaced
  by $w$.

  In summary, each vertex in $D$ is either already contained in $C$ or the
  result of application of Step (2).  Therefore, $D\coloneqq \DD{v}{G} \in
  \mathfrak{D}_G$ with $|D|>1$ if and only if $D$ can be obtained from
  $C\coloneqq \CC_{\lxt(G)}(v) \in \mathfrak{C}^*$ by applying Step (2).
  Hence, $\mathfrak{D}_G$ can be obtained from $\mathfrak{C}_{\lxt(G)}$ by
  first applying Step~(1) and then Step~(2).
  
  As shown above, each singleton $\{x\}\in \mathfrak{C}^*$ corresponds to
  the unique singleton $\{x\} \in \mathfrak{D}_G$ (in which case $x\in
  L(G)$) and each $\DD{v}{G} \in \mathfrak{D}_G$ with $|\DD{v}{G}|>1$ is
  obtained from $\CC_{\lxt(G)}(v) \in \mathfrak{C}^*$ by application of
  Step~(2). By Lemma~\ref{lem:desc-subset-iff-comparable}, $\DD{u}{G} \neq
  \DD{v}{G}$ for all distinct $u, v\in V(G)$ and, by
  Lemma~\ref{lem:G*-basic-properties}, $\CC_{\lxt(G)}(u)\neq
  \CC_{\lxt(G)}(v)$ for all distinct $u,v\in V(\lxt(G))$. Summarizing the
  latter arguments, there is a 1-to-1 correspondence between the elements
  in $\mathfrak{C}^*$ and the elements in $\mathfrak{D}_G$.
  
  We continue with showing that $(\mathfrak{D}_G,\subseteq)$ and
  $(\mathfrak{C}^*,\subseteq)$ are isomorphic. Based on the latter
  arguments, the map $\varphi \colon \mathfrak{D}_G \to \mathfrak{C}^*$
  where $\DD{v}{G} \in \mathfrak{D}_G$ is mapped to $\CC_{\lxt(G)}(v) \in
  \mathfrak{C}^*$ for all $v\in V(G)$ is bijective. Let $\DD{v}{G},
  \DD{w}{G}\in\mathfrak{D}_G$. By
  Lemma~\ref{lem:desc-subset-iff-comparable}, $\DD{v}{G} = \DD{w}{G}$ if
  and only if $v=w$. By Lemma~\ref{lem:G*-basic-properties}, $v=w$ if and
  only if $\CC_{\lxt(G)}(v)=\CC_{\lxt(G)}(w)$. Suppose now that $\DD{v}{G}
  \subsetneq \DD{w}{G}$ and thus, that $v\neq
  w$. Lemma~\ref{lem:desc-subset-iff-comparable} implies that $v\prec_G
  w$. By construction, $v,w\in V(\lxt(G))$ and $v\prec_{\lxt(G)} w$. This
  together with Lemma~\ref{lem:prec-subset} and
  \ref{lem:G*-basic-properties} implies that
  $\CC_{\lxt(G)}(v)\subsetneq\CC_{\lxt(G)}(w)$. Conversely, suppose that
  $\CC_{\lxt(G)}(v)\subsetneq\CC_{\lxt(G)}(w)$ and thus that $v\neq w$. By
  Lemma~\ref{lem:G*-basic-properties}, $\lxt(G)$ satisfies (PCC). However,
  $w\preceq_G v$ cannot holds because in this case
  Lemma~\ref{lem:prec-subset} would imply that
  $\CC_{\lxt(G)}(w)\subseteq\CC_{\lxt(G)}(v)$. Hence, $v\prec_G w$ must
  hold.  Lemma~\ref{lem:desc-subset-iff-comparable} implies that $\DD{v}{G}
  \subsetneq \DD{w}{G}$. In summary, $(\mathfrak{C}^*,\subseteq)$ and
  $(\mathfrak{D}_G,\subseteq)$ are isomorphic. Since, by
  Lemma~\ref{lem:desc-subset-iff-comparable}, $(\mathfrak{D}_G,\subseteq)$
  and $(V(G),\preceq_G)$ are isomorphic, we conclude that
  $(\mathfrak{C}^*,\subseteq)$ and $(V(G),\preceq_G)$ are isomorphic.
\end{proof}	 

A similar result can be obtained to show that $\mathfrak{C}_{\lxt(G)}$ can
be reconstructed from $\desccl{G}$ by applying the inverse of Step~(2), and
then the inverse of Step~(1).

\begin{lemma}\label{lem:dec-lxl-closed}
  Let $G$ be a DAG. If $\mathfrak{C}_{\lxt(G)}$ is closed, then
  $\desccl{G}$ is closed.
\end{lemma}
\begin{proof}
  Let $G$ be a DAG and assume that $\mathfrak{C}_{\lxt(G)}$ is closed. We
  denote with $\mathfrak{C}^*$ the set that is obtained from
  $\mathfrak{C}_{\lxt(G)}$ after application of Step~(1) in
  Lemma~\ref{lem:D-and-C_N*}, i.e., after removal of all singletons
  $\{x_v\}$ from $\mathfrak{C}_{\lxt(G)}$, where $x_v\in
  \overline{L}\coloneqq L(\lxt(G))\setminus L(G)$. Let $A,B\in \desccl{G}$
  be such that $A\cap B\neq \emptyset$. By Lemma~\ref{lem:D-and-C_N*},
  there is a 1-to-1 correspondence between the elements in $\desccl{G}$ and
  $\mathfrak{C}^*$. In particular, there is a unique element $A' \in
  \mathfrak{C}^*$ from which $A$ has been obtained after application of
  Step (2) on $A'$. In other words, all $x_v$ in $\overline L\cap A'$ have
  been replaced in $A'$ by the respective vertex $v\in V(G)$ to obtain
  $A$. Similarly, there is a unique $B' \in \mathfrak{C}^*$ from which $B$
  has been derived after application of Step (2) on $B'$. In particular,
  any $x_v \in A'\cap B'$ with $x_v \in \overline L$ has been replaced by
  $v$ in $A$ and $B$ while all other vertices $w\in A'\cap B'$ with $w
  \notin \overline L$ remain in $A$ and $B$ and therefore $C'\coloneqq
  A'\cap B'\neq\emptyset$.

  Since $A',B'\in \mathfrak{C}^*\subseteq \mathfrak{C}_{\lxt(G)}$ and
  $\mathfrak{C}_{\lxt(G)}$ is closed, it follows that $C'\in
  \mathfrak{C}_{\lxt(G)}$. Assume, for contradiction, $C' \notin
  \mathfrak{C}^*$. In this case, $C'= \{x_v\}$ for some $x_v \in \overline
  L$. Since $A',B'\in \mathfrak{C}^*$, it must hold that $A' \neq \{x_v\}$
  and $B'\neq \{x_v\}$. Consequently, the vertex $v$ to which $x_v$ was
  attached is an inner vertex in $G$, and there is a leaf $z\neq x_v$ with
  $z\in L(G)\subseteq L(\lxt(G))$ such that $z\prec_{\lxt(G)} v$. Moreover,
  there are vertices $u,w$ in $\lxt(G)$ for which $\CC_{\lxt(G)}(u) = A'$
  and $\CC_{\lxt(G)}(w) =B'$. Since $x_v\in A',B'$ and since $x_v$ has as
  its unique parent $v$, it follows that $v\preceq_{\lxt(G)} u,w$ and thus,
  that $z\prec_{\lxt(G)} u,w$. Therefore, $z\in A'$ and $z\in B'$ and, in
  particular, $z\in C'$. Since $z\neq x_v$ it follows that $C'\neq
  \{x_v\}$; a contradiction. Consequently, $C' \in \mathfrak{C}^*$ must
  hold. Due to the 1-to-1 correspondence between the elements in
  $\desccl{G}$ and $\mathfrak{C}^*$ we have an element $C\in\desccl{G}$
  that is obtained from $C'$ by replacing the vertices in $C'$ according to
  Step~(2) in Lemma~\ref{lem:D-and-C_N*}. This in particular implies that
  $C=A\cap B$. Since $C\in\desccl{G}$, the set $\desccl{G}$ is closed.
\end{proof}

Let us now turn to sets $\ANC_G(A)$ of all common ancestors of $A$ in $G$.
In the following we denote with $\mathfrak{A}_G$ the set comprising
$\ANC_G(\{v\})$ for all $v\in V(G)$ in a DAG $G$. By slight abuse of
notation we write $\ANC_G(v)$ instead of $\ANC_G(\{v\})$ when $\ANC_G$ is
applied to singleton sets.
\begin{observation}\label{obs:predec-G-descend-G<-}
  Let $G$ be a DAG. Then, for all vertices $v\in V$ it holds that
  $\DD{v}{G}\cap \ANC_G(v) = \{v\}$.  Moreover, every descendant (resp.,
  ancestor) of $u$ in $G$ is an ancestor (resp. descendant) of $u$ in
  $\rev(G)$ and, therefore, $\predecessorscl{G}=\desccl{\rev(G)}$ and
  $\desccl{G}=\predecessorscl{\rev(G)}$.
\end{observation}

\begin{lemma}\label{lem:weakly-closed-D-minimal-predecessor}
	Let $G = (V,E)$ be a DAG such that $\desccl{G}$ is closed. Then, for all $x,y\in V(G)$ with
	$\ANC_G(\{x,y\})\neq\emptyset$ it holds that $|\LCA_G(\{x,y\})|=1$.
\end{lemma}
\begin{proof}
	Let $G = (V,E)$ be a DAG such that $\desccl{G}$ is  closed and let $x,y\in V$ 
	be such that $\ANC\coloneqq \ANC_G(\{x,y\})\neq\emptyset$. Note that any common ancestor of
	$x$ and $y$ must be contained in $\ANC$. Since $\ANC\neq \emptyset$, we have $\LCA_G(\{x,y\})\neq
	\emptyset$. Assume, for contradiction, that $|\LCA_G(\{x,y\})|>1$ and let $w,w'\in \LCA_G(\{x,y\})$ be distinct. 
	By definition of LCAs, the elements in $\LCA_G(\{x,y\})$ must be
	pairwise $\preceq_G$-incomparable. By definition, $x,y \in D\coloneqq \DD{w}{G}\cap\DD{w'}{G}$.
	Since $\desccl{G}$ is closed, $D\in \desccl{G}$ and, hence, there is a vertex $v\in V$ such
	that $\DD{v}{G} = D$. Since, by definition, $v\in \DD{v}{G}$, it follows that $v$ is a descendant
	of $w$ and $w'$ and thus, $v\preceq_G w,w'$. This together with the fact that all elements in
	$\LCA_G(\{x,y\})$ are pairwise $\preceq_G$-incomparable implies that $v\notin \LCA_G(\{x,y\})$ and 
	thus, $v	\prec_G w,w'$; a contradiction to $w,w'\in \LCA_G(\{x,y\})$. Hence, $|\LCA_G(\{x,y\})|=1$ must hold.
\end{proof}

\begin{lemma}\label{lem:D-P-closed-equivalence}
	Let $G$ be a DAG. Then, $\desccl{G}$ is  closed if and only if $\predecessorscl{G}$ is  closed.
\end{lemma}
\begin{proof}
  Let $G$ be a DAG. Assume first that $\desccl{G}$ is closed. Let $A,B\in
  \predecessorscl{G}$ be such that $A\cap B\neq \emptyset$. Note that, by
  definition, there are vertices $u,v\in V(G)$ with with $A = \ANC_G(u)$
  and $B = \ANC_G(v)$.  In particular, $\ANC_G(\{u,v\}) = A\cap B \neq
  \emptyset$.  Hence, we can apply
  Lemma~\ref{lem:weakly-closed-D-minimal-predecessor} to conclude that
  $\LCA_G(\{u,v\})=\{w\}$ for some $w\in V(G)$.  By definition, $w\preceq_G z$
  for all $z\in \ANC_G(\{u,v\}) $ which implies that $\ANC_G(\{u,v\})
  \subseteq \ANC_G(w)$.  Moreover, $u,v\preceq_G w$ implies that every
  $z\in\ANC_G(w)$ satisfies $u,v\preceq_G z$ and, therefore, $z\in
  \ANC_G(\{u,v\}) $.  Hence, $\ANC_G(w) \subseteq \ANC_G(\{u,v\})$.  In
  summary, $\ANC_G(\{u,v\}) = A\cap B = \ANC_G(w) \in \predecessorscl{G}$.
  Thus, $\predecessorscl{G}$ is closed.
	
  For the converse, assume that $\predecessorscl{G}$ is closed. By
  Observation~\ref{obs:predec-G-descend-G<-},
  $\predecessorscl{G}=\desccl{\rev(G)}$. Hence, $ \desccl{\rev(G)}$ is
  closed. Since $\rev(G)$ is a DAG, we can re-use the previous arguments to
  conclude that $\predecessorscl{\rev(G)}$ is closed.  Since, by
  Observation~\ref{obs:predec-G-descend-G<-}, $\predecessorscl{\rev(G)} =
  \desccl{G}$ it follows that $\desccl{G}$ is closed.
\end{proof}

\begin{definition}\label{def:intermediary}
  For a DAG $G$ and two of its vertices $u,v\in V(G)$ let
  $\BB_G(u,v)\coloneqq \DD{u}{G} \cap \ANC_G(v) $
  be set of all intermediary
  vertices located between $u$ and $v$, i.e., the set of all vertices $w\in
  V(G)$ satisfying $v\preceq_G w\preceq_G u$.  The intermediary system is
  then $\allpathcl{G}\coloneqq\{\BB_G(u,v)\mid u,v\in V(G)\}$.
\end{definition}
	
It is easy to see that $\BB_G(u,v)$ contains all vertices that lie on any
directed $uv$-path. The map $(u,v)\mapsto \BB_G(u,v)$ coincides with the
\emph{directed all-path transit function} of the DAG, see
\cite{anil2024directed} for more details.  Note that, in a DAG $G$, it must
hold that $\BB_G(v,v) = \{v\}$ for all $v\in V(G)$. Hence, $\{v\} \in \allpathcl{G}$ 
for all $v\in V(G)$. However,  $\allpathcl{G}$ is in general not grounded
since, for any two vertices $u$ and $v$ for which
$u\prec_G v$ holds or where $u$ is $\preceq_G$-incomparable to $v$, we have
$\BB_G(u,v)=\emptyset$ and therefore, $\emptyset \in \allpathcl{G}$. 
Hence, $\allpathcl{G}$ is grounded if and only
if it contains a single vertex only. 
 Moreover, $\allpathcl{G}$ is not
necessarily rooted.

\begin{lemma}\label{lem:D-A-network-weakly-closed-equivalence}
  Let $N$ be a network, then $\predecessorscl{N}$ is closed if and only if
  $\allpathcl{N}$ is closed.
\end{lemma}
\begin{proof}
  Let $N$ be a network. Assume first that $\predecessorscl{N}$ is closed.
  By Lemma~\ref{lem:D-P-closed-equivalence}, $\desccl{N}$ is closed.
  Suppose there are vertices $x,y,u,v\in V(N)$ such that
  $\allpath{x}{y}{N}\cap\allpath{u}{v}{N}\neq\emptyset$.  We obtain
  \begin{align*}
    \emptyset \neq \allpath{x}{y}{N}\cap\allpath{u}{v}{N}
    =\ \DD{x}{N} \cap \ANC_N(y)  \cap \DD{u}{N} \cap \ANC_N(v)
    =\ &\underbrace{\DD{x}{N}\cap \DD{u}{N}}_{\neq\emptyset} \cap
     \underbrace{\ANC_N(y)  \cap \ANC_N(v)}_{\neq\emptyset} \\
    =\ & \DD{w}{N} \cap \ANC_N(z)
    = \allpath{w}{z}{N} \in \allpathcl{N},
  \end{align*}
  where we have used that $\predecessorscl{N}$ and $\desccl{N}$ are closed
  to justify the existence of vertices $w,z\in V(G)$ that satisfy
  $\DD{w}{N}=\DD{x}{N}\cap \DD{u}{N}$ and $\ANC_N(z)=\ANC_N(y) \cap
  \ANC_N(v)$.  Hence, if $\predecessorscl{N}$ is closed then
  $\allpathcl{N}$ is closed.
	
  Now suppose that $\allpathcl{N}$ is closed. Note that for every $v\in
  V(N)$, we have $\ANC_N(v)=\allpath{\rho}{v}{N}$, where $\rho$ denotes the
  unique root of $N$. Hence, for every $x,y\in V(N)$ such that
  $\ANC_N(x)\cap \ANC_N(y)\neq\emptyset$, it holds that
  $\allpath{\rho}{x}{N}\cap\allpath{\rho}{y}{N}\neq\emptyset$. Thus
  closedness of $\allpathcl{N}$ implies the existence of $w,z\in V(N)$ such
  that $\allpath{w}{z}{N}=\allpath{\rho}{x}{N}\cap\allpath{\rho}{y}{N}$. In
  particular, $\rho\in\allpath{w}{z}{N}$ holds. Since $\allpath{w}{z}{N} =
  \DD{w}{N}\cap \ANC_N(z)$, we have $\rho\in \DD{w}{N}$. This and the fact
  that $\rho$ is the unique $\preceq_N$-maximal vertex in $N$ implies that
  $w = \rho$, that is,
  $\allpath{w}{z}{N}=\allpath{\rho}{z}{N}=\ANC_N(z)$. Summarizing, we have
  \begin{equation*}
    \ANC_N(x)\cap
    \ANC_N(y)=\allpath{\rho}{x}{N}\cap\allpath{\rho}{y}{N}=
    \allpath{\rho}{z}{N}=\ANC_N(z)\in\predecessorscl{N},
  \end{equation*}
  and $\predecessorscl{N}$ is closed.
\end{proof}
	
\begin{theorem}\label{thm:global-lca-descendants}
  The following statements are equivalent for every network $N$.
  \begin{enumerate}[noitemsep,label=(\arabic*)]
  \item $N$ has the global lca-property.
  \item $\desccl{N}$ is   closed.
  \item $\predecessorscl{N}$ is  closed.
  \item $\allpathcl{N}$ is  closed.
  \end{enumerate}			
\end{theorem}
\begin{proof}
If $N$ has the global lca-property, then
Theorem~\ref{thm:char-global-lca-lxt} implies that $\mathfrak{C}_{\lxt(N)}$
is a closed clustering system.  By Lemma~\ref{lem:dec-lxl-closed},
$\desccl{N}$ is closed.  Suppose that $\desccl{N}$ is closed. Since $N$ is
a network, $\ANC_N(\{x,y\})\neq\emptyset$ for all $x,y\in V(N)$.  This and
Lemma \ref{lem:weakly-closed-D-minimal-predecessor} implies that
$|\LCA_G(x,y)|=1$ for all $x,y\in V(N)$. Hence, $N$ has the {global
  pairwise-lca-property and thus, by Theorem~\ref{thm:global-k2}, the}
global lca-property. In summary (1) and (2) are equivalent.  By
Lemma~\ref{lem:D-P-closed-equivalence}, (2) and (3) are equivalent and by
Lemma~\ref{lem:D-A-network-weakly-closed-equivalence}, (3) and (4) are
equivalent.
\end{proof}

\section{The Canonical Hasse-diagram of a DAG}
\label{sec:canonHasse}

{ A common theme throughout this contribution is
  that the existence of unique least common ancestors is closely tied to
  order-theoretic structures. For example, networks $N$ with global
  lca-property are precisely those DAGs whose underlying reachability poset
  $(V(N),\preceq_N)$ forms a join-semilattice
  (Theorem~\ref{thm:NgLCA-join.stuff}). Moreover, several of our
  characterizations are expressed in terms of Hasse diagrams of set systems
  associated with a DAG. From this perspective, it is natural to ask to
  what extent a DAG can be recovered from such Hasse-diagram
  representations.

  Intuitively, the Hasse diagrams of set systems derived from a DAG $G$
  should reflect the partial order $(V(G),\preceq_G)$ that governs ancestor
  relations and, hence, the structure of least common ancestors. This is
  particularly true for the poset $(\mathfrak{D}_G,\subseteq)$, which is
  isomorphic to $(V(G),\preceq_G)$ by
  Lemma~\ref{lem:desc-subset-iff-comparable}. As an immediate consequence
  of Lemmas~\ref{lem:desc-subset-iff-comparable} and
  \ref{lem:sf(G)=hasseDG}, we therefore have
  $\sf(G)\simeq\hasse(\mathfrak{D}_G,\subseteq)$. This relationship extends
  to the Hasse diagram of the poset $(\mathfrak{C}^*,\subseteq)$ obtained
  from $\mathfrak{C}_{\lxt(G)}$, due to its isomorphism with
  $(\mathfrak{D}_G,\subseteq)$ established in
  Lemma~\ref{lem:D-and-C_N*}. In contrast, the set system $\mathfrak{C}_G$,
  although directly defined in terms of descendant leaves and hence closely
  related to LCAs, often does not provide sufficient information to
  reconstruct $G$ (cf., e.g., Figure~\ref{fig:exmpl-lca-closed}).  This
  situation improves substantially after passing to the leaf-extended DAG
  $\lxt(G)$: by Theorem~\ref{thm:char-global-lca-lxt}, the LCA-properties
  of $G$ are faithfully reflected in $\lxt(G)$, while
  Lemma~\ref{lem:D-and-C_N*} shows that the cluster system
  $\mathfrak{C}_{\lxt(G)}$ encodes the descendant system $\mathfrak{D}_G$
  and thus, the poset $(\mathfrak{D}_G,\subseteq)$, which is isomorphic to
  $(V(G),\preceq_G)$.  Hence, the poset
  $(\mathfrak{C}_{\lxt(G)},\subseteq)$ yields a canonical order-theoretic
  representation of the reachability structure underlying the
  LCA structure of $G$.

  The aim of this section is to make this relationship explicit. More
  precisely, we show that, although $\mathfrak{C}_G$ alone is generally
  insufficient to recover $G$, the Hasse diagram of the poset
  $(\mathfrak{C}_{\lxt(G)},\subseteq)$ contains enough information to
  reconstruct $G$ up to shortcuts after a simple canonical pruning
  operation. Thus, this section complements the preceding LCA
  characterizations by showing how the order structure underlying least
  common ancestors can be recovered from an associated set of clusters.}

\begin{definition}
  A \emph{$G$-lopped} DAG $H$ is obtained form a tree-leaf-child DAG $G$ by
  removal of exactly one leaf-child $x_v$ of $v$ that has in-degree one,
  for all inner vertices $v\in V(H)$.  We denote by $\lop(G)$ an arbitrary
  $G$-lopped DAG.
\end{definition}

\begin{lemma}\label{lem:sflop=lopsf}
  Let $G$ be a tree-leaf-child DAG.  If $H$ and $H'$ are two $G$-lopped
  DAGs, then $H\simeq H'$. Moreover, $\sf(\lop(G)) \simeq \lop(\sf(G))$.
\end{lemma}
\begin{proof}
  Let $G = (V,E)$ be a {tree-leaf-child} DAG on $X$ and suppose that
  $H$ and $H'$ are two $G$-lopped DAGs. If $H=H'$, there is nothing to
  show. Assume that $H\neq H'$. In this case, $H$ is obtained from $G$ by
  removing a tree-leaf $x_v\in X$ and $H'$ is obtained from $G$ by removing
  a tree-leaf $x'_v\in X$ with $x_v\neq x_v'$ for at least one inner vertex
  $v$ of $G$. We call such particular inner vertices $v$
  \emph{$(H,H')$-leaf-deviating}. Note that all $(H,H')$-leaf-deviating
  vertices $v$ satisfy $v\in V(H)\cap V(H')$ as they must be non-leaf
  vertices. We can now define a bijective map $\varphi \colon V(H) \to
  V(H') $ by putting $\varphi(v) = v$ if $v\in V(H')$ and $\varphi(x_v) =
  x'_v$ for all $(H,H')$-leaf-deviating vertices $v\in V(H)$. In
  particular, the restriction of $\varphi$ to $V^0(G)\subseteq V(H)\cap
  V(H')$ is just the identity map. It is now an easy exercise to verify
  that $\varphi$ is an isomorphism between $H$ and $H'$.
  
  We continue with proving that $\sf(\lop(G)) \simeq \lop(\sf(G))$. Note
  that if there is a shortcut $e=(u,v)$ in $G$ for which the vertex $v$ is
  a leaf, then $v$ must be a hybrid in $G$ but may become a tree-leaf in
  $\sf(G)$ and thus, might be removed when considering $\lop(\sf(G))$. To
  address this issue, we mark, for each inner vertex $v$ in $G$, exactly
  one tree-leaf-child of $v$ red and call it $x_v$. The latter is possible
  since $G$ is tree-leaf-child. Let $S\subseteq E$ be the set of all
  shortcuts in $G$. The DAG $\sf(G)$ is obtained from $G$ by removal of all
  edges in $S$. Lemma~\ref{lem:shortcutfree} implies that $\sf(G)$ is
  $\preceq_G$-preserving{, i.e., $u\preceq_G v$ precisely if $u\preceq_{\sf(G)} v$
  for all $u,v\in V$.} Thus, $X = L(\sf(G))$ and $V^0 \coloneqq
  V^0(G) = V^0(\sf(G))$. Moreover, the in-degree of a vertex $w$ in
  $\sf(G)$ cannot be higher than the respective in-degree of $w$ in $G$,
  i.e., tree-leafs of $G$ remain tree-leafs in $\sf(G)$. Hence, all
  red-marked vertices in $G$ specify, in particular, a unique
  tree-leaf-child $x_v$ in $\sf(G)$ for all $v\in V^0$. We can thus obtain
  a particular $\sf(G)$-lopped DAG, denoted $\lop_{\text{\tiny
      red}}(\sf(G))$, by removing all red-marked vertices from
  $\sf(G)$. Similarly, we obtain a particular $G$-lopped DAG, denoted
  $\lop_{\text{\tiny red}}(G)$, by removing all red-marked vertices from
  $G$. It is easy to verify that the set of shortcuts in $\lop_{\text{\tiny
      red}}(G)$ coincides with $S$. Hence, to obtain $\sf(\lop_{\text{\tiny
      red}}(G))$ we can simply remove all edges in $S$. As we have applied
  identical operations to obtain $\sf(\lop_{\text{\tiny red}}(G))$ as well
  as $\lop_{\text{\tiny red}}(\sf(G))$ from $G$ it follows that
  $\sf(\lop_{\text{\tiny red}}(G)) = \lop_{\text{\tiny red}}(\sf(G))$. Now,
  as argued in the previous paragraph, all $G$-lopped DAGs are
  isomorphic. Hence, $\lop(G)\simeq \lop_{\text{\tiny red}}(G)$ for any
  $\lop(G)$. In addition,
  $\sf(H)$ is uniquely determined for any DAG $H$ which together with
  $\lop(G)\simeq \lop_{\text{\tiny red}}(G)$ implies that
  $\sf(\lop(G))\simeq\sf(\lop_{\text{\tiny red}}(G))$. By analogous
  arguments, it holds that $\lop(\sf(G))\simeq\lop_{\text{\tiny
      red}}(\sf(G))$ for every possible $\lop(\sf(G))$. Taking the latter
  two arguments together with $\sf(\lop_{\text{\tiny red}}(G)) =
  \lop_{\text{\tiny red}}(\sf(G))$ we can conclude that $\sf(\lop(G))
  \simeq \lop(\sf(G))$ always holds.
\end{proof}

{We close this section by providing Proposition~\ref{prop:N*lcarel},
  which shows that the order-theoretic structure underlying ancestor
  relations---and hence the structure from which least common ancestors
  arise---can be canonically recovered from the cluster system of the
  leaf-extended DAG $\lxt(G)$, even when the original cluster system
  $\mathfrak{C}_G$ is insufficient for this purpose.}

\begin{proposition}\label{prop:N*lcarel}
  Let $G$ be a DAG and $\lxt(G)$ its leaf-extended version. Then, $\sf(G)
  \simeq \lop(\hasse(\mathfrak{C}_{\lxt(G)},\subseteq))$. In other words,
  every DAG $G$ is isomorphic to a graph obtained from 
  $\lop(\hasse(\mathfrak{C}_{\lxt(G)},\subseteq))$ by adding a possibly
    empty set of shortcuts.
\end{proposition}
\begin{proof}
  Set $\hasse\coloneqq \hasse(\mathfrak{C}_{\lxt(G)},\subseteq)$. By
  Lemma~\ref{lem:G*-basic-properties}, $\sf(\lxt(G))$ is regular. Thus,
  $\sf(\lxt(G)) \simeq \hasse(\mathfrak{C}_{\sf(\lxt(G))},\subseteq)$. 
  {Note that $\sf(\lxt(G))$ is uniquely determined and that, 
	by  Lemma~\ref{lem:shortcutfree},}
  $\mathfrak{C}_{\sf(\lxt(G))} = \mathfrak{C}_{\lxt(G)}$. Consequently,
  $\sf(\lxt(G))\simeq \hasse$. By Lemma~\ref{lem:sflop=lopsf}, $G \simeq
  \lop(\lxt(G))$.  Hence, $\sf(G) \simeq \sf(\lop(\lxt(G))) \simeq
  \lop(\sf(\lxt(G))) \simeq \lop(\hasse)$, and $G$ is thus isomorphic to a
  DAG obtained from $\lop(\hasse)$ by adding some {$k\geq 0$} shortcuts.
\end{proof}

\section{Concluding Remarks} 
\label{sec:outlook} 
In this contribution, we have characterized global lca-networks or,
equivalently, global lca-DAGs, i.e., those DAGs that admit a unique LCA for
every subset of their vertices. It does not come as a big surprise that
global lca-DAGs $G$ are precisely those DAGs where the underlying poset
$(V(G),\preceq_G)$ is a join-semilattice. However, as shown here, global
lca-DAGs (and therefore join-semilattices) exhibit additional rich
structural properties: they are characterized by certain topological
minors, can be constructed stepwise from a single-vertex graph, and are
precisely those DAGs where the LCA is uniquely defined for each pair of
vertices. In addition, global lca-DAGs can be recognized in polynomial time
and are exactly those networks $N$ for which their descendant-system
$\desccl{N}$ (resp., ancestor-system $\predecessorscl{N}$ or
intermediary system $\allpathcl{N}$) is closed. Finally, we have shown that
each DAG $G$ can be uniquely reconstructed from the Hasse diagram of the
clustering system $\mathfrak{C}_{\lxt(G)}$ of its leaf-extended version
$\lxt(G)$.

The characterization of global lca-DAGs in terms of \holju graphs, i.e.,
the stepwise construction of such DAGs from a single-vertex graph, is of
particular interest. As shown in Lemma~\ref{lem:network-construction},
every network $N = (V, E)$ can be constructed from a single-vertex graph by
stepwise addition of a new vertex $x$, which is connected to \emph{some}
subset $W \subseteq V$. The characterization of global lca-DAGs then
involves restricting this set $W$ in each of the construction steps. It is
of interest to understand in more detail which DAGs with specified
properties admit a similar construction. In particular, how does the
addition of $x$, linked to the vertices in $W$, influence properties such
as the existence of LCAs or clustering systems with desired properties? For
instance, the clustering system of a level-1 network has been characterized
as those set systems that are closed and satisfy a property called (L)
\cite{Hellmuth:22q}. In particular, level-1 networks are \holju graphs
(cf.\ Cor~\ref{cor:level-1=holju}). This raises the question of which
properties of $W$ can be relaxed to ensure only that $\mathfrak{C}_{G}$ is
closed. Similarly, what conditions must be imposed on $W$ so that the
resulting network is an lca-network but not necessarily a global lca-network?

In contrast to clustering systems, the closedness of the underlying
descendant system $\desccl{N}$, ancestor system $\predecessorscl{N}$, and
intermediary system $\allpathcl{N}$ characterizes the global lca-property. 
Clearly, these three set systems contain all the information
about the underlying network needed to reconstruct it (up to its potential
shortcuts) and thus are nearly as informative as the underlying adjacency
matrix. In practical applications, one is often faced with the problem that
only partial information about the underlying set systems is provided. For
example, consider a (shortcut-free) DAG $G$ and its descendant system
$\desccl{G}$. Now consider a set system $\mathfrak{F}$ such that each $F
\in \mathfrak{F}$ satisfies $F \subseteq D$ for some $D \in \desccl{G}$,
denoted by $\mathfrak{F} \subseteq \desccl{G}$.  In other words,
$\mathfrak{F}$ is obtained from $\desccl{G}$ by removing some elements from
certain $D \in \desccl{G}$. In this case, $\mathfrak{F}$ contains only
partial information about $\desccl{G}$. This raises the question: to what
extent can $G$ be recovered from $\mathfrak{F}$? While this question might
be difficult to answer in general, it may be feasible for DAGs equipped
with structural constraints (e.g., lca-networks or global lca-networks). The
latter question is also is closely related to the question of how much
information about $G$ is contained in the ``containment-closure''
$\mathrm{CCL}(\mathfrak{F}) = \cap_{H, \mathfrak{F} \subseteq \desccl{H}}
\desccl{H}$, which is the intersection of all $\desccl{H}$ over all DAGs
$H$ for which $\mathfrak{F} \subseteq \desccl{H}$ holds. A similar concept
has been studied for so-called triples in phylogenetic trees
\cite{GRUNEWALD2007521}.

The semi-lattice convexity is known to be a convex geometry where
every convex set has the extreme-point property \cite{vandeVel:93}, i.e.,
each convex set $K$ contains a point $k\in K$ such that $K\setminus\{k\}$
is again convex. Every convex set is therefore the convex hull of its
extreme points. In our case, the extreme points of the join-semilattice
$(D,\sqsubseteq)$ are precisely the least common ancestors of the
join-subsemilattices $(A,\sqsubseteq),\ A\subseteq D$ that are again
convex.  The co-points of an element $p$ of this convex geometry are
maximal convex join-subsemilattices that do not contain the element
$p$. It is of practical interest for network reconstruction to study to
what extend $(D,\sqsubseteq)$ can be reconstructed from co-points of
certain restricted number of leaves.

\subsection*{Funding}

All authors declare that they have no conflicts of interest.

This work was supported in part by
the DST, Govt.\ of India (Grant No.\ DST/INT/DAAD/P-03/ 2023 (G)), the
DAAD, Germany (Grant No.\ 57683501), and the CSIR-HRDG for the Senior
Research Fellowship (09/0102(12336)/2021-EMR-I).
PFS acknowledges the financial support by the Federal Ministry of Education
and Research of Germany (BMBF) through DAAD project 57616814 (SECAI, School
of Embedded Composite AI), and jointly with the S{\"a}chsische
Staatsministerium f{\"u}r Wissenschaft, Kultur und Tourismus in the
programme Center of Excellence for AI-research \emph{Center for Scalable
Data Analytics and Artificial Intelligence Dresden/Leipzig}, project
identification number: SCADS24B.

\bibliographystyle{spbasic}
\bibliography{lca-minor}

\end{document}